\documentclass[11pt,a4paper]{amsart}

\newtheorem{theo}{Theorem}[section]
\newtheorem{defi}{Definition}[section]
\newtheorem{coro}[theo]{Corollary}
\newtheorem{prop}[theo]{Proposition}
\newtheorem{lemma}[theo]{Lemma}
\newcommand{\C}{\mathbb{C}}

\newcommand{\CP}{\mathbb{CP}}
\newcommand{\dbar}{\bar{\partial}}
\newcommand{\J}{\mathcal{J}}

\begin{document}
\title[Estimated transversality in symplectic geometry]
{Estimated transversality in symplectic geometry and projective maps}
\author{Denis Auroux}
\address{Centre de Math\'ematiques, Ecole Polytechnique, 
F-91128 Palaiseau, France}
\email{auroux@math.polytechnique.fr}

\maketitle
\section{Introduction}
Since Donaldson's original work \cite{D1}, approximately holomorphic
techniques have proven themselves most useful in symplectic geometry and
topology, and various classical constructions from algebraic geometry have 
been extended to the case of symplectic manifolds \cite{A2,A3,D2,mps}.
All these results rely on an estimated 
transversality statement for approximately holomorphic sections of very 
positive bundles, obtained by Donaldson \cite{D1,D2}.
However, the arguments require transversality not only for sections
but also for their covariant derivatives, which makes it necessary to
painstakingly imitate the arguments underlying Thom's classical {\it strong 
transversality theorem} for jets.

It is our aim in this paper to formulate and prove a general result of
estimated transversality with respect to finite stratifications in jet 
bundles. The transversality properties obtained in the various
above-mentioned papers then follow as direct corollaries of this result,
thus allowing some of the arguments to be greatly simplified.
The result can be formulated as follows (see \S 2 and \S 3 for
definitions)~:

\begin{theo}
Let $(E_k)_{k\gg 0}$ be an asymptotically very ample sequence of locally
splittable complex vector bundles over
a compact almost-complex manifold $(X,J)$.
Let $\mathcal{S}_k$ be asymptotically holomorphic finite
Whitney quasi-stratifications of the holomorphic jet bundles $\J^r E_k$.
Finally, let $\delta>0$ be a fixed constant. Then there
exist constants $K$ and $\eta$ such that, given any
asymptotically holomorphic sections $s_k$ of $E_k$ over $X$,
there exist asymptotically holomorphic 
sections $\sigma_k$ of $E_k$ with 
the following properties for all $k\ge K$~:

$(1)$ $|\sigma_k-s_k|_{C^{r+1},g_k}<\delta$~;

$(2)$ the jet $j^r\sigma_k$ of $\sigma_k$ is $\eta$-transverse to 
the quasi-stratification $\mathcal{S}_k$.
\end{theo}

We start by introducing in \S 2 a general notion of ampleness over an
almost-complex manifold. Then, in \S 3 we define the notion of approximately
holomorphic quasi-stratification of a jet bundle. Theorem 1.1 and its 
one-parameter version are proved in \S 4. Finally, we discuss applications 
in \S 5.

\section{Ample bundles over almost-complex manifolds}

The most general setup in which one can try to define a notion of ampleness
is the following. Let $X$ be a compact $2n$-dimensional
manifold (possibly with boundary), endowed with an 
almost-complex structure $J$. In order to make estimates, we also endow
$X$ with a Riemannian metric $g$ compatible with $J$ (i.e. $J$ is 
$g$-antisymmetric).

\begin{defi}
Given positive constants $c$ and $\delta$,
a complex line bundle $L$ over $X$ endowed with a Hermitian metric
and a connection $\nabla^L$ is {\em $(c,\delta)$-ample} if its curvature 
2-form $F_L$ satisfies the inequalities $iF_L(v,Jv)\ge c\,g(v,v)$
for every tangent vector $v\in TX$, and $\sup |F_L^{0,2}|\le
\delta$.

A sequence of complex line bundles $L_k$ with metrics and connections is
{\em asymptotically very ample} if there exist fixed constants
$\delta$ and $(C_r)_{r\ge 0}$, and a sequence $c_k\to+\infty$, 
such that the curvature $F_k$ of $L_k$ satisfies the following properties~: 
$(1)$~$iF_k(v,Jv)\ge c_k\,g(v,v)$ for every tangent vector $v\in TX$;
$(2)$~$\sup |F_k^{0,2}|\le \delta\,c_k^{1/2}$; 
$(3)$~$\sup |\nabla^r F_k|\le C_r\,c_k$ $\forall r\ge 0$.
\end{defi}

Most of this definition is a natural extension to the almost-complex setup
of the classical notion of ampleness on a complex manifold. Because the
notion of holomorphic bundle is not relevant in the case of a non-integrable
complex structure, one should allow the curvature to contain a non-trivial 
$(0,2)$-part. However, because $F_L^{0,2}$ is an obstruction to the
existence of holomorphic sections, we need uniform bounds on this quantity
in order to hope for the existence of approximately holomorphic sections.

The last condition in the definition seems less natural and should largely
be considered as a technical assumption needed to obtain some control over 
the behavior of sections; it is likely that a suitable argument, possibly
involving plurisubharmonic techniques, could allow the bounds to be 
significantly weakened.%
\medskip

Observe that the curvature of a $(c,\delta)$-ample line bundle over $X$ 
defines, after multiplication by $\frac{i}{2\pi}$, a $J$-tame symplectic 
structure on $X$ with integral cohomology class ($J$ is compatible with 
this symplectic structure if and only if the curvature is of type $(1,1)$).

Conversely, assume that $X$ carries a $J$-compatible symplectic form $\omega$ 
with integral cohomology class, and choose $g$ to be the Riemannian metric
induced by $J$ and $\omega$. Then there exists a line bundle $L$ with first 
Chern class $c_1(L)=[\omega]$ and a connection with curvature $-2\pi i\omega$
on $L$. By construction the line bundle $L$ is $(2\pi,0)$-ample;
moreover, the line bundles $L^{\otimes k}$ with the induced connections are 
$(2\pi k,0)$-ample and define an asymptotically very ample sequence of line 
bundles. This example is by far the most interesting one for
applications, but many other situations can be considered as well.%
\medskip

In the rest of this section, we consider an asymptotically very ample 
sequence of line bundles over $X$, and study the properties of $L_k$ for 
large values of $k$. In order to make the estimates below easier to 
understand, we rescale the metric by setting $g_k=c_k\,g$, which amounts to
dividing by $c_k^{r/2}$ the norm of all $r$-tensors; the Levi-Civita 
connection is not affected by this rescaling. The bounds of Definition 2.1 
imply that~: $iF_k(v,Jv)\ge g_k(v,v)$; $|F_k^{0,2}|_{g_k}=O(c_k^{-1/2})$;
$|F_k|_{g_k}=O(1)$~; $|\nabla^rF_k|_{g_k}=O(c_k^{-1/2})$ $\forall r\ge 1$.
Also observe that $|\nabla^r J|_{g_k}=O(c_k^{-1/2})$ $\forall r\ge 1$.
Better bounds on higher-order derivatives are trivially available but we 
won't need them.

\begin{lemma}
Let $L_k$ be a sequence of asymptotically very ample line bundles $L_k$ over
$X$, and denote by $F_k$ the curvature of $L_k$. Let $\omega_k=iF_k$, and
let $c_k$ be the constants appearing in Definition 2.1. 
Denote by $\nabla$ the Levi-Civita connection associated to $g$.
Then, for large enough $k$ there exist $\omega_k$-compatible 
almost-complex structures $\tilde{J}_k$ such that 
$|\nabla^r(\tilde{J}_k-J)|_{g_k}=O(c_k^{-1/2})$ $\forall r\ge 0$.
\end{lemma}%

\begin{proof}
We construct $\tilde{J}_k$ locally; patching together the local 
constructions in order to obtain a globally defined almost-complex structure
still satisfying the same type of bounds is an easy task left to the reader
(recall that the space of $\omega_k$-compatible almost-complex structures is
pointwise contractible).

Let $e_1$ be a local tangent vector field of unit $g_k$-length
and with $|\nabla^r e_1|_{g_k}=O(c_k^{-1/2})$ $\forall r\ge 1$ (observe
that, because of the rescaling process, $(X,g_k)$ is almost flat for large 
$k$). We define $e'_1=Je_1$, and observe that $e'_1$ has unit $g_k$-length 
($J$ is $g$-unitary and hence $g_k$-unitary) and $\omega_k(e_1,e'_1)\ge 1$. 
Next, we proceed inductively, assuming that we have defined local vector 
fields $e_1,e'_1,\dots,e_m,e'_m$ with the following properties for all 
$i,j\le m$~: $e_1,\dots,e_m$ have unit $g_k$-length;
$\omega_k(e_i,e_j)=\omega_k(e'_i,e'_j)=0$; $\omega_k(e_i,e'_j)=0$ if $i\neq
j$; $\omega_k(e_i,e'_i)\ge 1$; $e'_i-Je_i\in\mathrm{span}(e_1,e'_1,\dots,
e_{i-1},e'_{i-1})$; $|e'_i-Je_i|_{g_k}=O(c_k^{-1/2})$; $|\nabla^r e_i|_{g_k}=
O(c_k^{-1/2})$ and $|\nabla^r e'_i|_{g_k}=O(c_k^{-1/2})$ $\forall r\ge 1$.

We choose $e_{m+1}$ to be a $g_k$-unit vector field which is 
$\omega_k$-orthogonal to $e_1,e'_1,\dots,e_m,e'_m$. The bound on $|\nabla^r
\omega_k|_{g_k}$ implies that we can choose $e_{m+1}$ in such a way 
that $|\nabla^r e_{m+1}|_{g_k}=O(c_k^{-1/2})$. Next, we define
$$e'_{m+1}=Je_{m+1}+\sum_{i=1}^m \frac{\omega_k(e'_i,Je_{m+1})e_i-
\omega_k(e_i,Je_{m+1})e'_i}{\omega_k(e_i,e'_i)}.$$
By construction, $\omega_k(e_i,e'_{m+1})=\omega_k(e'_i,e'_{m+1})=0$ for all
$i\le m$. Moreover, $\omega_k(e_{m+1},e'_{m+1})=\omega_k(e_{m+1},Je_{m+1})
\ge 1$.

Since $\omega_k(Je_i,e_{m+1})=0$, and because
$\omega_k(Je_i,e_{m+1})-\omega_k(e_i,Je_{m+1})$ is a component of
$\omega_k^{0,2}$, we have
$\omega_k(e_i,Je_{m+1})=O(c_k^{-1/2})$.
Similarly,
we have
$\omega_k(e'_i,Je_{m+1})=\omega_k(Je_i,Je_{m+1})+
\omega_k(e'_i-Je_i,Je_{m+1})$; the first term differs from
$\omega_k(e_i,e_{m+1})=0$ by a $(0,2)$-term and is therefore bounded
by $O(c_k^{-1/2})$, while the bound on $e'_i-Je_i$ implies that
the second term is also bounded by $O(c_k^{-1/2})$. Therefore we have
$\omega_k(e'_i,Je_{m+1})=O(c_k^{-1/2})$.
Finally, using the lower bound on $\omega_k(e_i,e'_i)$ we obtain that
$|e'_{m+1}-Je_{m+1}|_{g_k}=O(c_k^{-1/2})$.
Finally, it is trivial that $|\nabla^r e'_{m+1}|_{g_k}=O(c_k^{-1/2})$~;
therefore we can proceed with the induction process.

We now define the almost-complex structure $\tilde{J}_k$ by the
identities $\tilde{J}_k(e_i)=e'_i$ and $\tilde{J}_k(e'_i)=-e_i$. 
By construction, $\tilde{J}_k$ is compatible with $\omega_k$, and the 
corresponding Riemannian metric $\tilde{g}_k$ admits $e_1,e'_1,\dots,e_n,e'_n$
as an orthonormal frame. The required bounds on $\tilde{J}_k$ immediately
follow from the available estimates.
\end{proof}

Lemma 2.1 makes it possible to recover the main ingredients of Donaldson
theory in the more general setting described here. We now introduce some
basic definitions and results, imitating Donaldson's original work and
subsequent papers \cite{D1,A2}.

In what follows, $L_k$ is an asymptotically very ample sequence of line 
bundles over $X$, $c_k$ are the same constants as in Definition 2.1, and 
$g_k=c_k\,g$.

\begin{lemma}
Near any point $x\in X$, and for any value of $k$, there exist local complex
Darboux coordinates $(z_k^1,\dots,z_k^n):(X,x)\to(C^n,0)$ for the symplectic
structure $\omega_k=iF_k$, such that, denoting by $\psi_k$ the inverse of
the coordinate map, the following bounds hold uniformly in $x$ and $k$ over
a ball of fixed $g$-radius around $x$~:
$|z_k^i(y)|=O(dist_{g_k}(x,y))$;
$|\nabla^r\psi_k|_{g_k}=O(1)$ $\forall r\ge 1$; and, with respect to the 
almost-complex structure $J$ on $X$
and the canonical complex structure on $\C^n$,
$|\dbar\psi_k(z)|_{g_k}=O(c_k^{-1/2}+c_k^{-1/2}|z|)$, and
$|\nabla^r\dbar\psi_k(z)|_{g_k}=O(c_k^{-1/2})$ for all $r\ge 1$.
\end{lemma}%

\begin{proof}
The argument is very similar to that used by Donaldson \cite{D1}, except 
that one needs to be slightly more careful in showing that the various bounds
hold uniformly in $k$. 
Fix a point $x\in X$~: then we can find a neighborhood $U$ of $x$ and a local
coordinate map $\phi:U\to\C^n$, such that $U$ contains a ball of fixed 
uniform $g$-radius around $x$, and such that the expressions of $g$ and
$J$ in these local coordinates satisfy uniform bounds independently of $x$
(these uniformity properties follow from the compactness of $X$). A linear
transformation can be used to ensure that the differential of $\phi$ at the 
origin is $\C$-linear with respect to $J$. Next, we rescale the
coordinates by $c_k^{1/2}$ to obtain a new coordinate map $\phi_k:U\to\C^n$, 
in which $J$ coincides with the standard almost-complex structure at
the origin and has derivatives bounded by $O(c_k^{-1/2})$, while the
expression of $g_k$ is bounded between fixed constants and has derivatives
bounded by $O(c_k^{-1/2})$. 

Next, we observe that the bound on $|\omega_k|$ and the lower bound on
$\omega_k(v,Jv)$ imply that the expression of $\omega_k^{(1,1)}$
at the origin of the coordinate chart is bounded from above and below by 
uniform constants. Therefore, after composing $\phi_k$ with a suitable 
element of $GL(n,\C)$, we can assume without affecting the bounds on $J$ and 
$g_k$ that $(\phi_k^{-1})^*(\omega_k^{(1,1)})$ 
coincides with the standard K\"ahler form $\omega_0$ of $\C^n$ at the
origin. 

Define over $\phi_k(U)\subset\C^n$ the symplectic form
$\omega_1=(\phi_k^{-1})^*\omega_k$. By construction,
$\omega_1(0)-\omega_0(0)=O(c_k^{-1/2})$. Observe that,
in the chosen coordinates, the Levi-Civita connection of 
$g_k$ differs from the trivial connection by $O(c_k^{-1/2})$~; therefore,
the bounds on $|\nabla^r\omega_k|_{g_k}$ imply that the derivatives of 
$\omega_1$ are also bounded by $O(c_k^{-1/2})$, and that
$|\omega_1(z)-\omega_0(z)|=O(c_k^{-1/2}+c_k^{-1/2}|z|)$.

In particular, decreasing the size of $U$ by at most a fixed factor if
necessary, we obtain that the closed
$2$-forms $\omega_t=t\omega_1+(1-t)\omega_0$ over $\phi_k(U)$ are all
symplectic, and we can apply Moser's argument to construct in a controlled
way a symplectomorphism between a subset of $(\phi_k(U),\omega_1)$ and a 
subset of $(\C^n,\omega_0)$. More precisely, it follows immediately from
Poincar\'e's lemma that we can choose a $1$-form $\alpha$ 
such that $\omega_1-\omega_0=d\alpha$, and such that $\alpha(0)=0$,
$|\alpha(z)|=O(c_k^{-1/2}|z|+c_k^{-1/2}|z|^2)$,
$|\nabla\alpha(z)|=O(c_k^{-1/2}+c_k^{-1/2}|z|)$ and
$|\nabla^r\alpha(z)|=O(c_k^{-1/2})$ $\forall r\ge 2$.
Next, we define vector fields $X_t$ by the identity $i_{X_t}\omega_t=\alpha$; 
clearly $X_t$ and its derivatives satisfy the same bounds as~$\alpha$.

Integrating the flow of the vector fields $X_t$ we obtain diffeomorphisms 
$\rho_t$, and it is a classical fact that the map 
$\tilde\phi_k=\rho_1\circ\phi_k$ is a local symplectomorphism between
$(X,\omega_k)$ and $(\C^n,\omega_0)$ and therefore defines Darboux 
coordinates. Because $|z|=O(c_k^{1/2})$ over a ball of fixed $g$-radius
around $x$, the vector fields $X_t$ satisfy a uniform bound of the type 
$|X_t(z)|\le \lambda |z|$ for some constant $\lambda$, so that
$|\rho_t(z)|\le e^{\lambda t}|z|$, and therefore $\tilde\phi_k$ is 
well-defined over a ball of fixed $g$-radius around $x$.
Moreover, the bounds $|\nabla(\rho_1-\mathrm{Id})|=O(c_k^{-1/2}+
c_k^{-1/2}|z|)$, obtained by integrating the bounds on $\nabla\alpha$, and 
$|\dbar(\phi_k^{-1})|=O(c_k^{-1/2}|z|)$, obtained from the bounds on the
expression on $J$ in the local coordinates, imply that
$|\dbar(\tilde\phi_k^{-1})|=O(c_k^{-1/2}+c_k^{-1/2}|z|)$. Similarly, the
bounds $|\nabla^{r+1} \rho_1|=O(c_k^{-1/2})$ and
$|\nabla^r\dbar(\phi_k^{-1})|=O(c_k^{-1/2})$ for all $r\ge 1$ imply
that $|\nabla^r\dbar(\tilde\phi_k^{-1})|=O(c_k^{-1/2})$.
This completes the proof of Lemma 2.2.
\end{proof}

\begin{defi}
A family of sections of $L_k$ is {\em asymptotically $J$-holomorphic} for
$k\to\infty$ if there exist constants $(C_r)_{r\ge 0}$ such that every
section $s\in\Gamma(L_k)$ in the family satisfies at every point of $X$ 
the bounds $|\nabla^r s|_{g_k}\le C_r$ and $|\nabla^r\dbar_J s|_{g_k}\le 
C_r\,c_k^{-1/2}$ for all $r\ge 0$, where $\dbar_J$ is the $(0,1)$-part
of the connection on $L_k$.

A family of sections of $L_k$ has {\em uniform Gaussian decay} properties
if there exist a constant $\lambda>0$ and polynomials $(P_r)_{r\ge 0}$ with 
the following property~: for every section $s$ of $L_k$ in the family, there 
exists a point $x\in X$ such that for all $y\in X$ and for all $r\ge 0$, 
$|\nabla^r s(y)|_{g_k}\le P_r(d_k(x,y))\exp(-\lambda\,d_k(x,y)^2)$, where 
$d(.,.)$ is the distance induced by $g_k$.
\end{defi}

\begin{lemma}
For all large enough values of $k$ and for every point $x\in X$, there
exists a section $s_{k,x}^\mathrm{ref}$ of $L_k$ with the following 
properties~: $(1)$ the family of sections 
$(s_{k,x}^\mathrm{ref})_{x\in X,k\gg 0}$ is asymptotically $J$-holomorphic;
$(2)$ the family $(s_{k,x}^\mathrm{ref})_{x\in X,k\gg 0}$ has uniform 
Gaussian decay properties, each section $s_{k,x}^\mathrm{ref}$ being
concentrated near the point $x$;
$(3)$ there exists a constant $\kappa>0$ independent of $x$ and $k$ such
that $|s_{k,x}^\mathrm{ref}|\ge \kappa$ at every point of the ball of
$g_k$-radius $1$ centered at $x$.
\end{lemma}%

\begin{proof}
The argument is a direct adaptation of the proof of Proposition 11 in 
Donaldson's paper \cite{D1}. Pick a value of $k$ and a point $x\in X$. 
We work in the approximately $J$-holomorphic Darboux coordinates given 
by Lemma 2.2, and use a trivialization of $L_k$ in which the connection 
$1$-form becomes $\frac{1}{4}\sum (z_j d\bar{z}_j-\bar{z}_j dz_j)$. 
Then, we define a local section of $L_k$ by 
$s(z)=\exp(-\frac{1}{4}|z|^2)$ and observe that $s$ is
holomorphic with respect to the standard complex structure of $\C^n$.
Multiplying $s$ by a cut-off function which equals $1$ over the ball of
radius $c_k^{1/6}$ around the origin, we obtain a globally defined section
of $L_k$; because of the estimates on the Darboux coordinates one easily 
checks that the families of sections constructed in this way are 
asymptotically holomorphic and have uniform Gaussian decay 
properties \cite{D1}.
\end{proof}

We are also interested in working with higher rank bundles. The definition
of ampleness becomes the following~:

\begin{defi}
A sequence of complex vector bundles $E_k$ with metrics and connections is
{\em asymptotically very ample} if there exist constants
$\delta$, $(C_r)_{r\ge 0}$, and $c_k\to+\infty$, 
such that the curvature $F_k$ of $E_k$ satisfies the following properties~: 
$(1)$ $\langle iF_k(v,Jv).u,u\rangle\ge c_k\,g(v,v)\,|u|^2$, 
$\forall v\in TX$, $\forall u\in E_k$;
$(2)$ $\sup |F_k^{0,2}|_g\le \delta_r\,c_k^{1/2}$; 
$(3)$ $\sup |\nabla^r F_k|_g\le C_r\,c_k$ $\forall r\ge 0$.

A sequence of asymptotically very ample complex vector bundles $E_k$ 
with metrics $|.|_k$
and connections $\nabla_k$ is {\em locally splittable} if, given any point 
$x\in X$, there exists over a ball of fixed $g$-radius around $x$ a 
decomposition of $E_k$ as a direct sum $L_{k,1}\oplus\dots \oplus L_{k,m}$ 
of line bundles, such that the following properties hold~:
$(1)$ the $|.|_k$-determinant of a local frame consisting of unit length
local sections of $L_{k,1},\dots,L_{k,m}$ is bounded from
below by a fixed constant independently of $x$ and $k$;
$(2)$ denoting by $\nabla_{k,i}$ the connection on $L_{k,i}$ obtained by 
projecting $\nabla_{k|L_{k,i}}$ to $L_{k,i}$, and by $\nabla'_k$
the direct sum of the $\nabla_{k,i}$, the $1$-form 
$\alpha_k=\nabla_k-\nabla'_k$ (the non-diagonal part of $\nabla_k$)
satisfies the uniform bounds $|\nabla^r\alpha_k|_g=O(c_k^{r/2})$ $\forall r\ge
0$ independently of $x$.
\end{defi}

For example, if $E$ is a fixed complex vector bundle and $L_k$ are
asymptotically very ample line bundles, then the vector bundles $E\otimes
L_k$ are locally splittable and asymptotically very ample; so are direct
sums of vector bundles of this type.

Observe that, if $E_k$ is an asymptotically very ample sequence of locally
splittable vector bundles, then near any given point $x\in X$ the summands 
$L_{k,1},\dots,L_{k,m}$ are asymptotically very ample line bundles. 
Therefore, by Lemma 2.3 they admit asymptotically holomorphic sections 
$s_{k,x,i}^\mathrm{ref}$ with uniform Gaussian decay away from $x$. 
Moreover, these sections, which define a local frame for $E_k$, are 
easily checked to be asymptotically $J$-holomorphic not only as 
sections of $L_{k,i}$ but also as sections of $E_k$. 

\section{Estimated transversality in jet bundles}

\subsection{Asymptotically holomorphic stratifications}

Throughout this section, we will denote by $F_k$ be a sequence of complex 
vector bundles over $X$, or more generally fiber bundles with almost-complex
manifolds as fibers. We also fix, in a manner compatible with the 
almost-complex structures $J^v$ of the fibers, metrics $g^v$ on the fibers 
of $F_k$ and a connection on $F_k$. Finally, we fix a sequence of constants
$c_k\to +\infty$.

The connection on $F_k$ induces a splitting $TF_k=T^vF_k\oplus T^hF_k$ between
horizontal and vertical tangent spaces~; this splitting makes it possible to
define a metric $\hat{g}_k$ and an almost-complex structure $\hat{J}_k$ on the
total space of $F_k$, obtained by orthogonal direct sum of $g^v$ and $J^v$
on $T^vF_k$ together with the pullbacks of $J$ and $g_k=c_k\,g$ on 
$T^hF_k\simeq \pi^*TX$.

We want to consider approximately holomorphic stratifications of the 
fibers of $F_k$, depending in an approximately holomorphic way on the point
in the base manifold $X$. For simplicity, we assume that the topological 
picture is the same in every fiber of $F_k$, i.e.\ we restrict ourselves to
stratifications which are everywhere transverse to the fibers. We will denote
the strata by $(S^a_k)_{a\in A_k}$~; we assume that the number of strata is
finite. Each $S^a_k$ is a possibly non-closed submanifold in $F_k$, whose 
closure is obtained by adding other lower dimensional strata~: writing 
$b\prec a$ iff $S^b_k$ is contained in $\overline{S^a_k}$, we have
$$\partial S^a_k\stackrel{\mathrm{def}}{=}
\overline{S^a_k}-S^a_k=\bigcup_{b\prec a}S^b_k.$$
We only consider 
{\em Whitney stratifications}~; in particular, transversality to a given
stratum $S^a_k$ implies transversality over a neighborhood of $S^a_k$ to
all the strata whose closure contains $S^a_k$, i.e.\ all the $S^b_k$ for
$b\succ a$. Also note that we discard any open strata, as they are
irrelevant for transversality purposes~; so each $S^a_k$ has codimension
at least $1$.

\begin{defi}
Let $(M,J)$ be an almost-complex manifold, with a Riemannian metric, and let
$s$ be a complex-valued function over $M$ or a section of an almost-complex 
bundle with metrics and connection. Given two constants $C$ and $c$,
we say that $s$ is {\em $C^2$-approximately holomorphic with bounds $(C,c)$},
or {\em $C^2$-AH($C,c$)}, if it satisfies the following estimates~:
$$|s|+|\nabla s|+|\nabla\nabla s|\le C,\qquad
|\dbar s|+|\nabla\dbar s|\le C\,c^{-1/2}.$$
Moreover, given constants $c_k\to+\infty$, we say that a sequence 
$(s_k)_{k\gg 0}$ of functions or sections is 
{\em $C^2$-asymptotically holomorphic}, or {\em $C^2$-AH}, if there exists
a fixed constant $C$ such that each section $s_k$ is
$C^2$-AH$(C,c_k)$.
\end{defi}

\begin{defi}
Let $F_k$ be a sequence of almost-complex bundles over $X$, endowed with
metrics and connections as above. For all values of $k$, let 
$(S^a_k)_{a\in A_k}$ be finite Whitney stratifications of $F_k$~; assume
that the total number of strata is bounded by a fixed constant
independently of $k$, and that all strata are transverse to the 
fibers of $F_k$.

We say that this sequence of stratifications is {\em asymptotically
holomorphic} if, given any bounded subset $U_k\subset F_k$,
and for every $\epsilon>0$, there exist positive constants
$C_\epsilon$ and $\rho_\epsilon$ depending only on $\epsilon$ and on the
size of the subset $U_k$ but not on $k$,
with the following property. For every point $x\in U_k$ lying 
in a certain stratum $S_k^a$ and at $\hat{g}_k$-distance greater than 
$\epsilon$ from $\partial S_k^a=\overline{S_k^a}-S_k^a$,
there exist complex-valued functions $f_1,\dots,f_p$
over the ball $B=B_{\hat{g}_k}(x,\rho_\epsilon)$ with the following
properties~:

$(1)$ a local equation of $S_k^a$ over $B$ is
$f_1=\dots=f_p=0$~;

$(2)$ $|df_1\wedge\dots\wedge df_p|_{\hat{g}_k}$ is bounded from below by 
$\rho_\epsilon$ at every point of $B$~;

$(3)$ the restrictions of $f_i$ to each fiber of $F_k$ near $x$ are
$C^2$-AH$(C_\epsilon,c_k)$~;

$(4)$ for any constant $\lambda>0$, and for any local section $s$
of $F_k$ which is $C^2$-AH$(\lambda,c_k)$ with respect to the metric $g_k$ on 
$X$ and which intersects non-trivially the ball $B$, the function $f_i\circ s$
is $C^2$-AH$(\lambda C_\epsilon,c_k)$~; moreover, given a local 
$C^2$-AH$(\lambda,c_k)$ section $\theta$ of $s^*T^vF_k$, the functions 
$df_i\circ \theta$ are $C^2$-AH$(\lambda C_\epsilon,c_k)$~;

$(5)$ at every point $y\in B$ belonging to a stratum $S_k^b$ such that
$S^a_k\subset\partial S^b_k$, the norm of the 
orthogonal projection onto 
the normal space $N_yS_k^b$ of any unit length vector $v\in T_yF_k$ such that 
$df_1(v)=\dots=df_p(v)=0$ is bounded by
$C_\epsilon\,\mathrm{dist}_{\hat{g}_k}(y,S_k^a)$.
\end{defi}

These conditions on the stratification can be reformulated more
geometrically as follows. 
First, the strata must be uniformly transverse to the fibers of $F_k$, 
i.e.\ one requires the minimum angle \cite{mps} between
$TS^a_k$ and $T^vF_k$ to be bounded from below. Second, the submanifolds
$S^a_k\subset F_k$ must be asymptotically $\hat{J}_k$-holomorphic, 
i.e.\ $\hat{J}_k(TS^a_k)$ and $TS^a_k$ lie within $O(c_k^{-1/2})$ of each
other. Third, the curvature of $S^a_k$ as a submanifold of $F_k$ must
be uniformly bounded. Finally, the quantity measuring the lack of
$\hat{J}_k$-holomorphicity of $S^a_k$ must similarly vary in a controlled way.

We finish this section by introducing the notion of estimated transversality
between a section and a stratification. Observe that, given any submanifold
$N\subset M$, we can define over a neighborhood of $N$ a ``parallel'' 
distribution $D_N\subset TM$ by parallel transport of $TN$ in the normal 
direction to $N$. Also recall that the {\it minimum angle} between two linear
subspaces $U$ and $V$ is defined as the minimum angle between a vector 
orthogonal to $U$ and a vector orthogonal to $V$ \cite{mps}. The minimum
angle between $U$ and $V$ is non-zero if and only if they are transverse to 
each other, and in that case it can also be defined as the minimum angle 
between non-zero vectors orthogonal to $U\cap V$ in $U$ and $V$.%

\begin{defi}
Given a constant $\eta>0$, we say that a section $s$ of a vector bundle
carrying a metric and a connection is $\eta$-transverse to $0$ if, 
at every point $x$ such that $|s(x)|\le \eta$, the covariant derivative
$\nabla s(x)$ is surjective and admits a right inverse of norm less than
$\eta^{-1}$.

Fix a constant $\eta>0$, and a section $s$ of a bundle carrying a
metric and a finite Whitney stratification $\mathcal{S}=(S^a)_{a\in A}$
everywhere transverse to the fibers. We say that $s$ is 
$\eta$-transverse to the stratification $\mathcal{S}$ if, at every point
where $s$ lies at distance less than $\eta$ from some stratum $S^a$, 
the graph of $s$ is transverse to the parallel
distribution $D_{S^a}$, with a minimum angle greater than $\eta$.

Finally, we say that a sequence of sections is uniformly tranverse to $0$
(resp.\ to a sequence of stratifications) if there exists a fixed constant
$\eta>0$ such that all sections in the sequence are $\eta$-transverse to
$0$ (resp.\ the stratifications).
\end{defi}

Note that the above condition of transversality of the section $s$ to each
stratum $S^a$ is only well-defined outside of a small neighborhood of
the lower-dimensional strata contained in $\partial S^a$~; however, near
these strata the assumption that $\mathcal{S}$ is Whitney makes
transversality to $S^a$ a direct consequence of the $\eta$-transversality to 
the lower-dimensional strata.

Another way in which uniform transversality to a stratification can be
formulated is to use local equations of the strata, as in Definition 3.2.
One can then define $\eta$-transversality as follows~:
at every point where $s$ lies at distance less than $\eta$ from $S^a$,
and considering local equations $f_1=\dots=f_p=0$ of $S^a$ such that
each $|df_i|$ is bounded by a fixed constant and 
$|df_1\wedge\dots\wedge df_p|$ is bounded from below by a fixed constant,
the function $(f_1\circ s,\dots,f_p\circ s)$ with values in $\C^p$ must be 
$\eta$-transverse to $0$.
The two definitions are equivalent up to changing the constant $\eta$
by at most a bounded factor.

\subsection{Quasi-stratifications in jet bundles}

Let $E_k$ be an asymptotically very ample sequence of locally splittable
rank $m$ vector bundles over the compact almost-complex manifold $(X,J)$. 
We can introduce the {\em holomorphic jet bundles}
$$\J^rE_k=\bigoplus_{j=0}^r 
\bigl(T^*X^{(1,0)}\bigr)^{\otimes j}_\mathrm{sym}\otimes E_k.$$

More precisely, the holomorphic part of the $r$-jet of a section $s$ of 
$E_k$ is defined inductively as follows~: $T^*X^{(1,0)}$ and $E_k$, as 
complex vector bundles carrying a connection over an almost-complex
manifold, are endowed with $\partial$ operators (the $(1,0)$ part of the 
connection)~; the $r$-jet of $s$ is $j^r s=(s,\partial_{E_k}
s,\partial_{T^*X^{(1,0)}\otimes E_k}(\partial_{E_k} s)_\mathrm{sym},\dots)$.

Observe that, because the almost-complex structure $J$ is not integrable and
because the curvature of $E_k$ is not necessarily of type $(1,1)$, 
the derivatives of order $\ge 2$ are not symmetric tensors,
but rather satisfy equality relations involving curvature terms and
lower-order derivatives. However, we will only consider the symmetric 
part of the jet~; for example, the $2$-tensor component of $j^r s$ is
defined by $(\partial\partial s)_\mathrm{sym}(u,v)=\frac{1}{2}
(\langle\partial(\partial s),u\otimes v\rangle+\langle\partial(\partial s),
v\otimes u\rangle)$. Note that, anyway, in
the case of asymptotically holomorphic sections, the antisymmetric terms
are bounded by $O(c_k^{-1/2})$, because the $(2,0)$ curvature terms and 
Nijenhuis tensor are bounded by $O(c_k^{-1/2})$.

The metrics and connections on $TX$ and on $E_k$ naturally induce 
Hermitian metrics and connections on $\J^rE_k$ (to define the metric we use 
the rescaled metric $g_k$ on $X$). In fact, it is easy to
see that the vector bundles $\J^rE_k$ are asymptotically very ample.

Recall that, near any given point $x\in X$, there exist
local approximately holomorphic coordinates~; besides a local 
identification of $X$ with $\C^n$, these coordinates also provide an
identification of $T^*X^{(1,0)}$ with ${T^*\C^n}^{(1,0)}$.
Moreover, by Lemma 2.3 there exist asymptotically 
holomorphic sections $s_{k,x,i}^\mathrm{ref}$ of $E_k$ with Gaussian decay 
away from $x$ and defining a local frame in $E_k$. Using these sections
to trivialize $E_k$, we can locally identify $\J^rE_k$ with a space of jets
of holomorphic $\C^m$-valued maps over $\C^n$. Observe however that,
when we consider the holomorphic parts of jets of approximately holomorphic
sections of $E_k$, the integrability conditions normally satisfied by jets
of holomorphic functions only hold in an approximate sense.

In general, the various possible choices of trivializations of $\J^rE_k$ 
differ by approximately holomorphic diffeomorphisms of $\C^n$ and also by 
the action of approximately holomorphic local sections of the automorphism 
bundle $\mathrm{GL}(E_k)$. However, when $E_k$ is of the form $\C^m\otimes
L_k$ where $L_k$ is a line bundle, the only automorphisms of $E_k$ which we
need to consider are multiplications by complex-valued functions.

Denote by $\J^r_{n,m}$ the space of $r$-jets of holomorphic maps from
$\C^n$ to $\C^m$~: pointwise, the identifications of the fibers
of $\J^r E_k$ with $\J^r_{n,m}$ given by local trivializations differ
from each other by the action of $GL_n(\C)\times GL_m(\C)$ (or $GL_n(\C)
\times\C^*$ when $E_k=\C^m\otimes L_k$), where $GL_n(\C)$ corresponds to
changes in the identification of $T^*X^{(1,0)}$ with ${T^*\C^n}^{(1,0)}$ and 
$GL_m(\C)$ or $\C^*$ corresponds to changes in the trivialization of $E_k$.
Some stratifications of $\J^r_{n,m}$ are invariant under the actions of 
$GL_n(\C)$ and $GL_m(\C)$ (resp.\ $\C^*$). Given such a
stratification it becomes easy to construct an asymptotically holomorphic
sequence of finite Whitney stratifications of $\J^rE_k$, modelled in each
fiber on the given stratification of $\J^r_{n,m}$. Many important 
examples of asymptotically holomorphic stratifications, and in a certain 
sense all the geometrically relevant ones, are obtained by this
construction (see Proposition 3.1 below).\medskip

We also wish to consider cases where the available structure is not exactly
a Whitney stratification but behaves in a similar manner with respect to
transversality. We call such a structure a ``Whitney quasi-stratification''.
Given a submanifold $S\subset\J^r_{n,m}$, one can introduce the
subset $\Theta_S$ of all points $\sigma\in S$ such that there exists a
holomorphic $(r+1)$-jet whose $r$-jet component is $\sigma$ and which,
considered as a $1$-jet of $r$-jets, intersects $S$ transversely at
$\sigma$. For example, if $S$ is the subset of all jets
$(\sigma_0,\dots,\sigma_r)$ such that $\sigma_0=0$, the subset $\Theta_S$
consists of those jets such that $\sigma_0=0$ and $\sigma_1$ is surjective.

Similarly, when $S$ is a submanifold in
$\J^rE_k$, we can view an element of $\J^{r+1}E_k$ as the 
holomorphic $1$-jet of a section of $\J^r E_k$. More precisely, for
any point $x\in X$, we can associate to any 
$\sigma=(\sigma_0,\dots,\sigma_{r+1})
\in(\J^{r+1}E_k)_x$ the $1$-jet at $x$ of a local section $\tilde\sigma$
of $\J^r E_k$, such that $\tilde\sigma(x)=(\sigma_0,\dots,\sigma_r)$,
$(\partial\tilde\sigma(x))^\mathrm{sym}=(\sigma_1,\dots,\sigma_{r+1})$,
$(\partial\tilde\sigma(x))^\mathrm{antisym}=0$, and
$\dbar\tilde\sigma(x)=0$ (in this definition, $\partial\tilde\sigma(x)
\in T^*X^{1,0}\otimes(\bigoplus (T^*X^{1,0})^{\otimes j}_\mathrm{sym}\otimes
E)$ is decomposed into a symmetric part and an antisymmetric part).
Then, we define $\Theta_S$ as the set of points of $S$ for which
there exists an element $\sigma\in\J^{r+1}E_k$ such that the
corresponding 1-jet $\tilde\sigma$ in $\J^rE_k$ intersects $S$ transversely
at the given point.
For example, if $S$ is the
set of $r$-jets $(\sigma_0,\dots,\sigma_r)$ such that $\sigma_0=0$, then
$\Theta_S$ is the set of $r$-jets such that $\sigma_0=0$
and $\sigma_1$ is surjective. Also observe that $\Theta_S$ is always empty
when the codimension of $S$ is greater than $n$.

\begin{defi}
Given a finite set $(A,\prec)$ carrying a binary relation without cycles 
(i.e., $a_1\prec\dots\prec a_p\Rightarrow a_p\not\prec a_1$),
a finite Whitney quasi-stratification of $\J^r_{n,m}$ indexed by $A$ is a 
collection $(S^a)_{a\in A}$ of smooth submanifolds of $\J^r_{m,n}$, not 
necessarily mutually disjoint, with the following
properties~: $(1)$ 
$\partial S^a=\overline{S^a}-S^a\subseteq \bigcup_{b\prec a} S^b$~;
$(2)$ given any point $p\in\partial S^a$, there exists $b\prec a$ such
that $p\in S^b$ and such that either $p\not\in\Theta_{S^b}$
or $S^b\subset\partial S^a$ and the Whitney regularity condition is 
satisfied at all points of $S^b$.
\end{defi}

Similarly, we can define the notion of asymptotically holomorphic finite
Whitney quasi-stratifications of $\J^r E_k$. This is similar to
Definition 3.2, except that the collections $(S^a_k)$ are
quasi-stratifications rather than stratifications, 
i.e.\ $\partial S^a_k\subseteq\bigcup_{b\prec a} S^b_k$, and for every
$p\in \partial S^a_k$ there exists $b\prec a$ such that either 
$p\in S^b_k-\Theta_{S^b_k}$ or $p\in S^b_k\subset\partial S^a_k$~; in the 
latter case the Whitney condition is required. Also, observe that condition 
$(5)$ in Definition 3.2 is only required in the second case,
and not for all $b$ such that $a\prec b$.

It is important to understand that the notion of quasi-stratification is 
merely an attempt at simplifying the framework for applications of 
Theorem 1.1. In fact, most quasi-stratifications can be refined into genuine
stratifications by suitably subdividing the strata into smaller pieces.
However, by definition these modifications occur at points of $\J^r E_k$
that no generic jet can hit, thus making them utterly irrelevant to
transversality.

\begin{prop}
Let $\mathcal{S}=(S^a)_{a\in A}$ be a finite Whitney quasi-stratification of
$\J^r_{n,m}$ by complex submanifolds, invariant under the action of
$GL_n(\C)\times GL_m(\C)$ or $GL_n(\C)\times \C^*$. 
Let $E_k$ be an asymptotically
very ample sequence of rank $m$ complex vector bundles over $X$, trivialized 
near every point by suitable choices of local asymptotically holomorphic 
coordinates and sections. Assume that
$\mathcal{S}_k=(S^a_k)_{a\in A}$ are quasi-stratifications of $\J^r E_k$ such
that, in each local trivialization, the intersection of $S^a_k$
with every fiber becomes identified with $S^a$. Then the sequence of
quasi-stratifications $\mathcal{S}_k$ is asymptotically holomorphic.
\end{prop}

The proof of this result is easy and left to the reader~; the
independence on $k$ of the model holomorphic quasi-stratification of 
$\J^r_{n,m}$ and the availability of asymptotically holomorphic local
trivializations of $\J^r E_k$ (Lemma 2.2 and Lemma 2.3) immediately yield
the necessary estimates on the strata of $\mathcal{S}_k$. The only important
point to observe is that, because the strata of $\mathcal{S}$ are
$GL_m(\C)$-invariant (resp.\ $\C^*$-invariant), the local trivializations
identifying $S^a_k$ with $S^a$ also identify $\Theta_{S^a_k}$ with
$\Theta_{S^a}$. This is e.g.\ due to the fact that, up to a suitable
change in the choice of the local coordinates on $X$ and local reference 
sections of $L_k$, i.e.\ up to a local gauge transformation, we can assume
that the connection on $\J^r E_k$ agrees at a given point $x\in X$ with the 
trivial connection on $\J^r_{n,m}$~; above $x$, the identification of 
$(r+1)$-jets with $1$-jets of $r$-jets then becomes the same in $\J^r E_k$
as in $\J^r_{n,m}$, so that the definitions of $\Theta_S$ in $\J^r_{n,m}$
and in $\J^r E_k$ agree with each other.

Various examples of applications of Proposition 3.1 will be given in \S 5.

Finally, we state a one-parameter version of Theorem 1.1. Consider a
continuous one-parameter family $(J_t)_{t\in [0,1]}$ of almost-complex
structures on $X$, and a one-parameter family of asymptotically holomorphic 
finite Whitney (quasi)-stratifications 
$(\mathcal{S}_{k,t})_{k\gg 0,t\in [0,1]}$ of
almost-complex bundles $F_{k,t}$ over $(X,J_t)$. We say that the 
(quasi)-stratifications $\mathcal{S}_{k,t}$ depend continuously on $t$ if the 
following one-parameter version of Definition 3.2 is true~:
for every $\epsilon>0$, there exist constants $\rho_\epsilon$ and
$C_\epsilon$ with the following property. Given any continuous path 
$(x_t)_{t\in [t_1,t_2]}$ of points all belonging to the fibers of $F_{k,t}$
above a same point in $X$, and assuming that all the points $x_t$ belong
to certain strata $S_{k,t}^a$ while lying at distance more than $\epsilon$
from $\partial S_{k,t}^a$, there exist for all $t\in [t_1,t_2]$ 
complex-valued functions $f_{1,t},\dots,f_{p,t}$ defined over the ball
$B_{\hat{g}_{k,t}}(x_t,\rho_\epsilon)$ and depending continuously on $t$,
satisfying the various properties of Definition 3.2 for all values of $t$.

With this understood, the result is the following~:

\begin{theo}
Let $(J_t)_{t\in [0,1]}$ be a continuous one-parameter family of 
almost-complex structures on the compact manifold $X$, and let 
$(E_{k,t})_{k\gg 0,t\in [0,1]}$ be a family of complex vector bundles over
$X$ endowed with metrics and connections depending continuously on $t$ and
such that the sequence $E_{k,t}$ is asymptotically 
very ample and locally splittable over $(X,J_t)$ for all $t$.
Let $\mathcal{S}_{k,t}$ be asymptotically holomorphic finite
Whitney quasi-stratifications of $\J^r E_{k,t}$
depending continuously on $t$.
Finally, let $\delta>0$ be a fixed constant. Then there
exist constants $K$ and $\eta$ such that, given any one-parameter family of
asymptotically holomorphic sections $s_{k,t}$ of $E_{k,t}$ over $X$
depending continuously on $t$, there exist asymptotically holomorphic 
sections $\sigma_{k,t}$ of $E_{k,t}$, depending continuously on $t$, with 
the following properties for all $k\ge K$ and for all $t\in [0,1]$~:

$(1)$ $|\sigma_{k,t}-s_{k,t}|_{C^{r+1},g_k}<\delta$~;

$(2)$ the jet $j^r\sigma_{k,t}$ of $\sigma_{k,t}$ is $\eta$-transverse to 
$\mathcal{S}_{k,t}$.
\end{theo}

\section{Proof of the main result}

The proof of Theorem 1.1 is quite similar to the arguments in previous
papers \cite{A2,A3,D1,D2,mps}. It relies heavily on the fact that 
the estimated transversality of the $r$-jet of a section to a given
submanifold of the jet bundle is a local and $C^{r+1}$-open property 
in the following sense \cite{A2}. Given a submanifold $S$ of $\J^r E_k$,
a constant $\eta>0$ and a point $x\in X$, say that a section $s$ of $E_k$
satisfies the property $\mathcal{P}(S,\eta,x)$ if either the $r$-jet
$j^r s(x)$ lies at distance more than $\eta$ from $S$ or 
$j^r s$ is $\eta$-transverse to $S$ at $x$ (in the sense of Definition 3.3,
i.e.\ the minimum angle between the graph of $j^r s$ and the parallel
distribution to $S$ is at least $\eta$). 
The property $\mathcal{P}(S,\eta,x)$ depends only on the $(r+1)$-jet of 
$s$ at $x$ (``locality''). Moreover, if $s$ satisfies 
$\mathcal{P}(S,\eta,x)$, then any section $\sigma$ such that
$|j^{r+1}\sigma(x)-j^{r+1}s(x)|<\epsilon$ satisfies
$\mathcal{P}(S,\eta-C\epsilon,x)$, where $C$ is some fixed 
constant involving only the curvature bounds of $S$
(``openness'').

A first consequence is that Theorem 1.1 can be proved by successively
perturbing the given sections $s_k$ in order to ensure transversality to
the various strata. To show this, we first remark that, given any index
$b$, the uniform transversality of $j^r s_k$ to all the strata $S^a_k$ 
with $a\prec b$ implies its uniform transversality to $S^b_k$ over a
neighborhood of $\partial S^b_k$.

Indeed, first consider a pair of indices $a\prec b$ such that
$S^a_k\subset \partial S^b_k$. By
condition $(5)$ of Definition 3.2, near a point of $S^a_k$ the tangent 
space to $S^b_k$ almost contains the parallel distribution to $TS^a_k$~;
therefore, there exists a constant $\kappa$ (independent of $a$ and $b$) 
such that, for any small $\alpha>0$, the $\alpha$-transversality of 
$j^r s_k$ to $S^a_k$ implies its $\frac{\alpha}{4}$-transversality to 
$S^b_k$ over the $\kappa\alpha$-neighborhood of $S^a_k$. 
Next, consider a pair of indices $a\prec b$ and a point 
$p\in\partial S^b_k\cap(S^a_k-\Theta_{S^a_k})$~: in this case,
if the graph of $j^r s_k$ is $\alpha$-transverse to $S^a_k$
but intersects the ball of radius $\frac{\alpha}{2}$ around $p$,
we can find an approximately holomorphic section 
$\sigma_k$ of $E_k$ differing from $s_k$ by less than 
$\frac{3\alpha}{4}$ and whose jet goes through $p$. By definition of
$\Theta_{S^a_k}$, all lifts of $p$ in $\J^{r+1}E_k$, including
$j^{r+1}\sigma_k$, correspond to local sections which intersect $S^a_k$ 
non-transversely; because the antisymmetric and antiholomorphic terms in
$\nabla(j^r\sigma_k)$ are smaller than $O(c_k^{-1/2})$, the minimum angle 
between $j^r\sigma_k$ and $S^a_k$ at $p$ is bounded by
$O(c_k^{-1/2})$. However, since $\sigma_k$ is close to $s_k$, its $r$-jet
should be $\frac{\alpha}{4}$-transverse to $S^a_k$, which gives a
contradiction. Therefore, $j^rs_k$ remains at distance more than
$\frac{\alpha}{2}$ from $p$; this implies the
$\frac{\alpha}{4}$-transversality to $S^b_k$ of $j^r s_k$ over the
$\frac{\alpha}{4}$-neighborhood of every point of 
$(S^a_k-\Theta_{S^a_k})\cap\partial S^b_k$. Since these are the only
two possible cases near the boundary of $S^b_k$, the uniform
transversality of $j^r s_k$ to $S^a_k$ for all $a\prec b$ implies its
uniform transversality to $S^b_k$ near $\partial S^b_k$.

Now, extend the binary relation $\prec$ on the
set of strata of each $\mathcal{S}_k$ into a total order relation $<$,
so that the indices in $A_k$ can be identified with integers and 
the closure of a given stratum consists only of strata appearing before it.
Assume that a first perturbation by 
less than $\delta_0=\frac{\delta}{2}$ makes it possible to obtain for large 
$k$ the $\eta_1$-transversality of $j^rs_k$ to the first stratum $S^1_k$, 
for some constant $\eta_1$ independent of $k$. Next, let $\delta_1$ be
a constant sufficiently smaller than $\delta$ and $\eta_1$ (but
independent of $k$), and assume that a perturbation
by at most $\delta_1$ allows us to obtain the
$\eta_2$-transversality of $j^r s_k$ to the second stratum $S^2_k$ outside of
the $\frac{1}{4}\kappa\eta_1$-neighborhood of $\partial S^2_k$, 
for some constant $\eta_2$. Because this new perturbation is small enough, 
the resulting sections remain $\frac{\eta_1}{2}$-transverse to $S^1_k$~; 
also, by the above observation this automatically implies the estimated
transversality to $S^2_k$ of $j^r s_k$ near the points of $\partial S^2_k
\subseteq S^1_k$.

We can continue in this way until all strata have been considered~; each
perturbation added to ensure estimated transversality to a new stratum 
outside of a small fixed size neighborhood of its boundary
is chosen small enough in order not to affect the previously
obtained transversality properties. 

The fact that estimated transversality is local and open also makes it
possible to reduce to a purely local setup, using a
globalization principle due to Donaldson \cite{D1} and which
can be formulated as follows (Proposition 3 of \cite{A2})~:

\begin{prop}
Let $\mathcal{P}_k(\eta,x)_{x\in X,\eta>0,k\gg 0}$ be local and 
$C^{r+1}$-open properties of sections of $E_k$ over $X$. 
Assume that there exist constants $c$, $c'$ and $\nu$ such that, given any 
$x\in X$, any small enough $\delta>0$, and asymptotically holomorphic 
sections $s_k$ of $E_k$, there exist, for all large enough $k$, 
asymptotically holomorphic sections $\tau_{k,x}$ of $E_k$ with the following 
properties~: {\rm (a)} $|\tau_{k,x}|_{C^{r+1},g_k}<\delta$,
{\rm (b)} the sections $\frac{1}{\delta}\tau_{k,x}$ have uniform Gaussian
decay away from $x$, and
{\rm (c)} the sections $s_k+\tau_{k,x}$ satisfy the property 
$\mathcal{P}_k(\eta,y)$
for all $y\in B_{g_k}(x,c)$, with $\eta=c'\delta\log(\delta^{-1})^{-\nu}$.

Then, given any $\alpha>0$ and asymptotically holomorphic sections $s_k$ of 
$E_k$, there exist, for all large enough $k$, asymptotically holomorphic 
sections $\sigma_k$ of $E_k$ such that $|s_k-\sigma_k|_{C^{r+1},g_k}<\alpha$ 
and the sections $\sigma_k$ satisfy $\mathcal{P}_k(\epsilon,x)$ 
$\forall x\in X$ for some $\epsilon>0$ independent of $k$.
\end{prop}

Proposition 4.1 is in fact slightly stronger than the previous results, 
as the notion of asymptotic holomorphicity has been extended to a more
general framework in \S 2, but the argument remains strictly the same.

With this result, we are reduced to the problem of finding a localized
perturbation of $s_k$ near a given point $x$ in order to ensure
transversality to a given stratum. More precisely, fix an index $a\in A_k$
in each stratification, and remember that, from the previous steps of the
inductive argument, we can restrict ourselves to considering only 
asymptotically holomorphic sections whose jet is $\gamma$-transverse to the 
strata $S^b_k$ for $b<a$, for some fixed constant $\gamma$ (this constant
$\gamma$ is half of the transversality estimate obtained in the previous 
step~; by assumption we only consider perturbations which are small
enough to preserve $\gamma$-transversality to the previous strata). 
With this understood, say that a section $s_k$ satisfies 
$\mathcal{P}_k(\eta,x)$ if either $j^r s_k(x)$ lies at distance more than
$\eta$ from $S^a_k$, or $j^r s_k(x)$ lies at distance less than
$\frac{1}{4}\kappa\gamma-\eta$ from $\partial S^a_k$,
or $j^r s_k$ is $\eta$-transverse to $S^a_k$ at $x$. We want to show that 
the assumptions of Proposition 4.1 are satisfied by these properties.

Fix a point $x\in X$ and a constant $0<\delta<\frac{1}{20}\kappa\gamma$,
and consider asymptotically holomorphic sections $s_k$ of $E_k$. 
First, if $j^rs_k(x)$ lies at 
distance less than $\frac{3}{20}\kappa\gamma$ from a point of
$\partial S^a_k\cap S_k^b$ for some $b\prec a$, 
then the uniform bounds on covariant derivatives of 
$s_k$ imply that the graph of $j^r s_k$ remains within distance less than 
$\frac{1}{5}\kappa\gamma$ of this point over a ball of fixed radius $c_1$ 
(independent of $k$, $x$ or $\delta$) around $x$. So, the property 
$\mathcal{P}_k(\frac{1}{20}\kappa\gamma,y)$ holds at every point
$y\in B_{g_k}(x,c_1)$, and no perturbation is needed. In the rest of the
argument, we can therefore assume that $j^r s_k(x)$ lies at distance at 
least $\frac{3}{20}\kappa\gamma$ from $\partial S^a_k$.

Let $\epsilon=\frac{1}{10}\kappa\gamma$, and let $\rho_\epsilon$ be the
radius appearing in Definition 3.2. Without loss of generality we can
assume that $\rho_\epsilon<\epsilon$. Assume that $j^r s_k(x)$ lies at
distance more than $\frac{1}{2}\rho_\epsilon$ from $S^a_k$. Then, the
bounds on covariant derivatives of $s_k$ imply that the graph of
$j^r s_k$ remains at distance more than $\frac{1}{4}\rho_\epsilon$ from
$S^a_k$ over a ball of fixed radius $c_2$ around $x$, and therefore that
$s_k$ satisfies $\mathcal{P}_k(\frac{1}{4}\rho_\epsilon,y)$ at every
point $y\in B_{g_k}(x,c_2)$. No perturbation is needed.

Therefore, we may assume that $j^r s_k(x)$ lies at distance less than
$\frac{1}{2}\rho_\epsilon$ from a certain point $u_0\in S^a_k$. We may also 
safely assume that $\delta<\frac{1}{4}\rho_\epsilon$. One easily
checks that $u_0$ lies at distance more than $\epsilon$ from $\partial
S^a_k$. So we can find complex-valued
functions $f_1,\dots,f_p$ over the ball $B_{\hat{g}_k}(u_0,\rho_\epsilon)$
such that a local equation of $S^a_k$ is $f_1=\dots=f_p=0$ and satisfying
the various properties listed in Definition 3.2. Let $c_3$ be a fixed
positive constant (independent of $k$, $x$ and $\delta$) such that the
graph of $j^r s_k$ over $B_{g_k}(x,c_3)$ is contained in $B_{\hat{g}_k}(u_0,
\frac{3}{4}\rho_\epsilon)$, and define the $\C^p$-valued function
$h=(f_1\circ j^r s_k,\dots,f_p\circ j^rs_k)$ over $B_{g_k}(x,c_3)$. 
By property $(4)$ of
Definition 3.2, the function $h$ is $C^2$-approximately holomorphic.

Recall from Lemma 2.2 that there exist local approximately holomorphic
$\omega_k$-Darboux coordinates $z_1,\dots,z_n$ over a neighborhood of $x$
in $X$. Also recall from Lemma 2.3 that there exist approximately holomorphic
sections $s_{k,x,i}^\mathrm{ref}$ of $E_k$ with Gaussian decay away from $x$ 
and defining a local frame in $E_k$. For any $(n+1)$-tuple 
$I=(i_0,i_1,\dots,i_n)$ with $1\le i_0\le m$, $i_1,\dots,i_n\ge 0$, 
and $i_1+\dots+i_n\le r$, we define $s_{k,x,I}^\mathrm{ref}=z_1^{i_1}
\dots z_n^{i_n} s_{k,x,i_0}^\mathrm{ref}$. Clearly, these 
sections of $E_k$ are asymptotically holomorphic and have uniform Gaussian 
decay away from $x$~; moreover it is easy to check that their $r$-jets
define a local frame in $\J^r E_k$ near $x$. After multiplication by
a suitable fixed constant factor, we can also assume that
$|s_{k,x,I}^\mathrm{ref}|_{C^{r+1},g_k}\le \frac{1}{p}$. 
For each tuple $I$, define a $\C^p$-valued function $\Theta_I$ by
$\Theta_I(y)=(df_1(j^r s_k(y)).j^r s_{k,x,I}^\mathrm{ref}(y),\dots,
df_p(j^r s_k(y)).j^r s_{k,x,I}^\mathrm{ref}(y))$. The functions $\Theta_I$
measure the variations of the function $h$ when small multiples of the
localized perturbations $s_{k,x,I}^\mathrm{ref}$ are added to $s_k$~; by
condition $(4)$ of Definition 3.2, they are $C^2$-asymptotically
holomorphic.

The fact that the jets of $s_{k,x,I}^\mathrm{ref}$ define a frame of
$\J^rE_k$ near $x$ implies, by condition $(2)$ of Definition 3.2, that
the values $\Theta_I(x)$ generate all of $\C^p$. Moreover, for $1\le i\le p$
there exist complex constants $\lambda_{I,i}$ with 
$\sum_I |\lambda_{I,i}|\le 1$ such that, defining the linear combinations 
$\sigma_{k,x,i}=\sum_{I}\lambda_{I,i} s_{k,x,I}^\mathrm{ref}$ and
$\Theta_i=\sum_I \lambda_{I,i} \Theta_I$, the quantity
$|\Theta_1(x)\wedge\dots\wedge\Theta_p(x)|$ is larger than some
fixed positive constant $\beta>0$ depending only on $\epsilon$ (and not
on $k$, $x$ or $\delta$). The uniform bounds on derivatives imply that,
for some fixed constant $0<c_4<c_3$, the norm of 
$\Theta_1\wedge\dots\wedge \Theta_p$ remains larger than 
$\frac{1}{2}\beta$ at every point of $B_{g_k}(x,c_4)$. Therefore, over
this ball we can express $h$ in the form $h=\mu_1\Theta_1+\dots+\mu_p
\Theta_p$, and the $\C^p$-valued function $\mu=(\mu_1,\dots,\mu_p)$ is
easily checked to be $C^2$-AH as well.

Finally, use once more the local approximately holomorphic coordinates to 
identify $B_{g_k}(x,c_4)$ with a neighborhood of the origin in $\C^n$. 
After rescaling the coordinates by a fixed constant factor, we can assume
that this neighborhood contains the ball $B^+$ of radius $\frac{11}{10}$
around the origin in $\C^n$, and that there exists a fixed constant
$0<c_5<c_4$ such that the inverse image of the unit ball $B$ in $\C^n$ 
contains $B_{g_k}(x,c_5)$. Composing $\mu$ with the coordinate map, we
obtain a $\C^p$-valued function $\tilde\mu$ over $B^+$~; by construction
$\tilde\mu$ is $C^2$-AH.

We may now use the following local result, due to Donaldson \cite{D2}
(the case $p=1$ is an earlier result of Donaldson \cite{D1}; the
comparatively much easier case $p>n$ is handled in \cite{A2})~:

\begin{prop}[Donaldson \cite{D2}]
Let $f$ be a function with values in $\C^p$ defined over the ball $B^+$ 
of radius $\frac{11}{10}$ in $\C^n$. Let $\delta$ be a constant with
$0<\delta<\frac{1}{2}$, and let $\eta=\delta\log(\delta^{-1})^{-\nu}$ where
$\nu$ is a suitable fixed integer depending only on $n$ and
$p$. Assume that $f$ satisfies the following bounds over $B^+$:
$$|f|\le 1,\qquad |\dbar f|\le \eta,\qquad |\nabla\dbar f|\le\eta.$$
Then, there exists $w\in\C^p$, with $|w|\le \delta$, such that $f-w$
is $\eta$-transverse to $0$ over the interior ball $B$ of radius $1$.
\end{prop}

Let $\eta=\delta\log(\delta^{-1})^{-\nu}$ as in the statement of
the proposition, and observe that, if $k$ is large enough, the 
antiholomorphic derives of $\tilde\mu$, which are bounded by a fixed 
multiple of $c_k^{-1/2}$, are smaller than $\eta$.
Therefore, if $k$ is large enough we can apply Proposition 4.2 (after a 
suitable rescaling to ensure that $\tilde\mu$ is bouded by $1$) and find a 
constant $w=(w_1,\dots,w_p)\in\C^p$,
smaller than $\delta$, such that $\tilde{\mu}-w$ is $\eta$-transverse to $0$
over the unit ball $B$. Going back through the coordinate map, this implies
that $\mu-w$ is $c'_1\eta$-transverse to $0$ over $B_{g_k}(x,c_5)$ for some 
fixed constant $c'_1$. Multiplying by the functions $\Theta_1,\dots,
\Theta_p$, we obtain that $h-(w_1\Theta_1+\dots+w_p\Theta_p)$ is
$c'_2\eta$-transverse to $0$ over $B_{g_k}(x,c_5)$ for some fixed constant
$c'_2$. 

Let $\tau_{k,x}=-(w_1 \sigma_{k,x,1}+\dots+
w_p \sigma_{k,x,p})$~: by construction, the sections $\tau_{k,x}$
of $E_k$ are asymptotically holomorphic, their norm is bounded by $\delta$, 
and they have uniform Gaussian decay properties. Let
$\tilde{s}_k=s_k+\tau_{k,x}$, and observe that by construction the graph
of $j^r\tilde{s}_k$ over $B_{g_k}(x,c_3)$ is contained in
$B_{\hat{g}_k}(u_0,\rho_\epsilon)$. Define $\tilde{h}=(f_1\circ
j^r\tilde{s}_k,\dots,f_p\circ j^r\tilde{s}_k)$~; by construction, and
because of the bounds on second derivatives of $f_1,\dots,f_p$, we have
the equality
$\tilde{h}=h-(w_1\Theta_1+\dots+w_p\Theta_p)+O(\delta^2)$.
If $\delta$ is assumed to be small enough, the quadratic term in this
expression is much smaller than $\eta$~; therefore, under this assumption
we get that $\tilde{h}$ is $c'_3\eta$-transverse to $0$ over
$B_{g_k}(x,c_5)$ for some fixed constant $c'_3$. Finally, recalling the
characterization of estimated transversality to a submanifold defined by
local equations given at the end of \S 3.1, we conclude that
the graph of $j^r\tilde{s}_k$ is $c'_4\eta$-transverse to $S^a_k$ over
$B_{g_k}(x,c_5)$ for some fixed constant $c'_4$, i.e.\ $\tilde{s}_k$
satisfies the property $\mathcal{P}_k(c'_4\delta\log(\delta^{-1})^{-\nu},y)$ 
at every point $y\in B_{g_k}(x,c_5)$.

Putting together the various possible cases (according to the distance
between $j^rs_k(x)$ and $S_k^a$ or its boundary), we obtain
that the properties $\mathcal{P}_k$ satisfy the assumptions of Proposition
4.1. Therefore, for all large values of $k$ a small perturbation can be
added to $s_k$ in order to achieve uniform transversality to $S_k^a$ away
from $\partial S_k^a$. The inductive argument described at the
beginning of this section then makes it possible to complete the proof
of Theorem 1.1.
\medskip

The proof of Theorem 3.2 follows the same argument, but for one-parameter
families of sections. One easily checks that the various results of \S 2
(Lemma 2.1, 2.2, 2.3) remain valid for families of
objects depending continuously on a parameter $t\in [0,1]$.
Moreover, Propositions 4.1 and 4.2 also extend to the one-parameter
case \cite{D2,A2}. So we only need to check that the argument used above
to verify that the properties $\mathcal{P}_k$ satisfy the assumptions of
Proposition 4.1 extends to the case of one-parameter families.

As before, fix a stratum $S^a_{k,t}$ in each stratification, a constant 
$\delta>0$, a point $x\in X$, and asymptotically holomorphic sections 
$s_{k,t}$ of $E_{k,t}$. With the same notations as above, let 
$\Omega_k\subset [0,1]$ be the set of values of $t$ such that 
$j^r s_{k,t}(x)$ lies at distance more than $\frac{3}{20}\kappa\gamma$ from
$\partial S^a_{k,t}$, and within distance
$\frac{1}{2}\rho_\epsilon$ from $S^a_{k,t}$. Let $\Omega_k^-\subset\Omega_k$
be the set of values of $t$ such that $j^r s_{k,t}(x)$ lies at distance more
than $\frac{1}{5}\kappa\gamma$ from $\partial S^a_{k,t}$ and within 
distance $\frac{1}{4}\rho_\epsilon$ from $S^a_{k,t}$.
Observe that, if $t\not\in\Omega_k^-$, a certain uniform transversality 
property with respect to $S^a_{k,t}$ is already satisfied by $j^r s_{k,t}$
over a small ball centered at $x$, and therefore no specific perturbation
is needed~: if $x$ lies within distance $\frac{1}{5}\kappa\gamma$ from
$\partial S^a_{k,t}$, then
$\mathcal{P}_k(\frac{1}{40}\kappa\gamma,y)$ is satisfied at every point
of a ball of fixed radius, while if $x$ lies at distance more than
$\frac{1}{4}\rho_\epsilon$ from $S^a_{k,t}$ then $\mathcal{P}_k(\frac{1}{8}
\rho_\epsilon,y)$ holds over a ball of fixed radius around $x$. 
Even better, if $\delta$ is small enough compared to $\gamma$ and 
$\rho_\epsilon$, then any perturbation of $s_{k,t}$ by less than $\delta$
still satisfies a similar transversality property (with decreased estimates).

For $t$ in $\Omega_k$, the proximity of $j^r s_{k,t}(x)$ to $S^a_{k,t}$
makes it possible to locally define complex-valued functions $f_{1,t},\dots,
f_{p,t}$ depending continuously on $t$ and such that a local equation of 
$S^a_{k,t}$ is $f_{1,t}=\dots=f_{p,t}=0$ (recall the definition of the 
continuous dependence of the stratifications $\mathcal{S}_{k,t}$ upon 
the parameter $t$ given in \S 3.2). This lets us define as above 
the function $h_t=(f_{1,t}\circ j^r s_{k,t},\dots,f_{p,t}\circ j^r 
s_{k,t})$, depending continuously on $t$. As in the non-parametric case, we
can construct asymptotically holomorphic sections $s_{k,x,t,I}^\mathrm{ref}$
of $E_{k,t}$, with Gaussian decay away from $x$ and defining local frames in 
$\J^r E_{k,t}$, simply by multiplying the sections of Lemma 2.3 by
polynomials of degree at most $r$ in the local coordinates (all these
sections depend continuously on $t$).
We can then find linear combinations $\sigma_{k,x,t,1},\dots,\sigma_{k,x,t,p}$
of the sections $s_{k,x,t,I}^\mathrm{ref}$, with constant coefficients
depending continuously on $t$, such that, denoting by $\Theta_{t,i}$ the 
$\C^p$-valued functions expressing the variations of 
$h_t$ upon adding small multiples of $\sigma_{k,x,t,i}$ to $s_{k,t}$,
the norm of $\Theta_{t,1}\wedge\dots\wedge\Theta_{t,p}$ is bounded from 
below at $x$ and over a small ball surrounding it. 

Constructing the functions $\tilde\mu_t$
as in the proof of Theorem 1.1 and applying the one-parameter version of
Proposition 4.2, we obtain, if $k$ is large enough, a continuous one-parameter
family of constants $w_t\in\C^p$, depending continuously on $t\in\Omega_k$ 
and bounded by $\delta$ for all $t$, such that $\tilde\mu_t-w_t$ is
$\eta$-transverse to $0$ over the unit ball in $\C^n$. It follows that,
denoting by $\tau_{k,x,t}$ the asymptotically holomorphic perturbations
$-(w_{t,1}\sigma_{k,x,t,1}+\dots+w_{t,p}\sigma_{k,x,t,p})$, bounded by
$\delta$, with Gaussian decay away from $x$, and depending continuously
on $t\in \Omega_k$, the sections $s_{k,t}+\tau_{k,x,t}$ satisfy the desired 
transversality property over a small ball centered at $x$. However these
perturbations are only well-defined for $t\in\Omega_k$. In order to extend
their definition to all values of $t$, let $\chi_k:[0,1]\to[0,1]$ be a
continuous cut-off function such that $\chi_k(t)=1$ for all $t\in \Omega_k^-$
and $\chi_k(t)=0$ for all $t\not\in \Omega_k$, and let $\tilde\tau_{k,x,t}=
\chi_k(t)\tau_{k,x,t}$ (for $t\not\in\Omega_k$ we set
$\tilde\tau_{k,x,t}=0$). For $t\in\Omega_k^-$ we have $\tilde\tau_{k,x,t}=
\tau_{k,x,t}$, so the sections $s_{k,t}+\tilde\tau_{k,x,t}$ satisfy the
required transversality property~; for $t\not\in\Omega_k^-$, the sections
$s_{k,t}$ already satisfy such a property and, because we have assumed
$\delta$ to be small enough, transversality is not affected by adding 
$\tilde\tau_{k,x,t}$. Therefore, the assumptions of Proposition 4.1 are
satisfied even in the one-parameter setting, and we can conclude the
argument in the same way as in the non-parametric case.%
\medskip

\noindent {\bf Remark.} In many cases, Theorems 1.1 and 3.2 can be proved 
without using Proposition 4.2 (estimated Sard lemma) in its full generality. 
Indeed, given suitable asymptotically holomorphic quasi-stratifications 
$\mathcal{S}_k$ of $\J^r E_k$, we can define quasi-stratifications 
$\tilde{\mathcal{S}}_k$ of $\J^{r+1}E_k$
in the following way. View each element of $\J^{r+1} E_k$ as a 1-jet in
$\J^r E_k$, as in \S 3.2; for each stratum $S^a_k$ of $\mathcal{S}_k$ with
codimension greater than $n$ in $\J^r E_k$, let $\tilde{S}^a_k$ be the
set of points in $\J^{r+1} E_k$ whose $r$-jet component belongs to $S^a_k$.
For each stratum $S^a_k$ of $\mathcal{S}_k$ with codimension $p\le n$, and
for each value $0\le i\le p-1$, let $\tilde{S}^{a,i}_k$ be the set of
points in $\J^{r+1} E_k$ whose $r$-jet component belongs to $S^a_k$ and
such that the corresponding element in $T^*X^{(1,0)}\otimes\J^r E_k$, after
projection to the normal space to $TS^a_k$, has rank equal to $i$. 
In other terms, the union of $\tilde{S}^{a,i}_k$ is the set of $(r+1)$-jets
which intersect $S^a_k$ non-transversely.

In a large number of examples, those of the $\tilde{S}^{a,i}_k$ which are
not empty are approximately holomorphic submanifolds of $\J^{r+1} E_k$,
transverse to the fibers and of codimension at least $n+1$. These
submanifolds determine finite Whitney quasi-stratifications $\tilde{\mathcal{S}}_k$ 
of $\J^{r+1} E_k$, satisfying properties similar to those of Definition 3.2
but with $C^1$ estimates only instead of $C^2$ bounds. Still, the same
argument as in the proof of Theorem 1.1 shows that, given
asymptotically holomorphic sections $s_k$ of $E_k$, small perturbations
can be added for large enough $k$ in order to ensure the uniform 
transversality of $j^{r+1}s_k$ to $\tilde{\mathcal{S}}_k$~; the argument
only uses Proposition 4.2 in the case $p>n$, where the proof becomes much
easier \cite{A2} and $C^1$ bounds are sufficient. Because all the strata are
of codimension greater than $n$, the $\eta$-transversality
of $j^{r+1}s_k$ to $\tilde{\mathcal{S}}_k$ simply means that the graph
of $j^{r+1}s_k$ remains at distance more than $\eta$ from the strata of 
$\tilde{\mathcal{S}}_k$. By definition of $\tilde{\mathcal{S}}_k$, this 
is equivalent to the uniform transversality of $j^r s_k$ to $\mathcal{S}_k$,
which was the desired result.

\section{Examples and applications}

We now consider various examples of (quasi)-stratifications to which we can 
apply Theorems 1.1 and 3.2. The fact that they are 
asymptotically holomorphic is in all cases a direct consequence of 
Proposition 3.1.

To make things more topological, we place ourselves in the case where 
the almost-complex structure $J$ on $X$ is tamed by a given symplectic form 
$\omega$. In this context, the various approximately $J$-holomorphic
submanifolds of $X$ appearing in the constructions are automatically
symplectic with respect to $\omega$. 
Moreover, remember that the space of $\omega$-tame or $\omega$-compatible
almost-complex structures on $X$ is contractible. In most
applications, asymptotically very ample bundles are 
constructed from line bundles with first Chern class proportional to
$[\omega]$~; in that situation, the ampleness properties of these bundles
do not depend on the choice of an $\omega$-compatible almost-complex 
structure $J$. Theorem 3.2 then implies that all the constructions described 
below are, for large enough values of $k$, canonical up to isotopy, 
independently of the choice of $J$. In the general case, the constructions
are still canonical up to isotopy, but the space of possible choices for
$J$ is constrained by the necessity for the bundles $E_k$ to be ample.

The first application is the construction of symplectic submanifolds
as zero sets of asymptotically holomorphic sections of vector bundles over
$X$, as initially obtained by Donaldson \cite{D1} and later extended to a
slightly more general setting \cite{A1}.

\begin{coro}
Let $(X,\omega)$ be a compact symplectic manifold endowed with an $\omega$-tame
almost-complex structure $J$, and let $E_k$ be an asymptotically very ample 
sequence of locally splittable vector bundles over $(X,J)$. Then, for all 
large enough values of $k$ there exist asymptotically holomorphic sections 
$s_k$ of $E_k$ which are uniformly transverse to $0$ and whose zero sets are 
smooth symplectic manifolds in $X$. Moreover these sections and submanifolds
are, for large $k$, canonical up to isotopy, indepedently of the chosen
almost-complex structure on $X$.
\end{coro}%

\begin{proof} 
Let $\mathcal{S}_k$ be the stratification of $\J^0E_k=E_k$ in which
the only stratum is the zero section of $E_k$ (these stratifications are
obviously asymptotically holomorphic). By Theorem 1.1, starting
from any asymptotically holomorphic sections of $E_k$ (e.g.\ the zero
sections) we can obtain for large $k$ asymptotically holomorphic sections
of $E_k$ which are uniformly transverse to $\mathcal{S}_k$, i.e.\ uniformly
transverse to $0$. It is then a simple observation that the zero sets of
these sections are, for large $k$, smooth approximately $J$-holomorphic
(and therefore symplectic) submanifolds of $X$ \cite{D1}. Finally, the
uniqueness of the construction up to isotopy is a direct consequence of
the one-parameter result Theorem 3.2 \cite{A1}.
\end{proof}

The next example is that of determinantal submanifolds as constructed by
Mu\~noz, Presas and Sols \cite{mps}.

\begin{coro}
Let $(X,\omega)$ be a compact symplectic manifold endowed with an $\omega$-tame
almost-complex structure $J$, let $L_k$ be an asymptotically very ample 
sequence of line bundles over $(X,J)$, and let $E$ and $F$ be complex
vector bundles over $X$. Then, for all large enough values of $k$ there exist
asymptotically holomorphic sections $s_k$ of $E^*\otimes F\otimes L_k$ such
that the determinantal loci
$\Sigma_i(s_k)=\{x\in X,\ \mathrm{rk}(s_k(x))=i\}$ are stratified symplectic 
submanifolds in $X$. Moreover these sections and submanifolds
are, for large $k$, canonical up to isotopy, indepedently of the chosen
almost-complex structure on $X$.
\end{coro}%

\begin{proof} 
Let $E_k=E^*\otimes F\otimes L_k$, and let $\mathcal{S}_k$ be the 
stratification of $\J^0E_k=E_k$ consisting of strata $S_k^i$,
$0\le i<\min(\mathrm{rk}\,E,\mathrm{rk}\,F)$, defined as follows~: viewing
the points of $E_k$ as elements of $\mathrm{Hom}(E,F)$ with coefficients in
$L_k$, each $S_k^i$ is the set of all elements in $E_k$ whose rank is equal to 
$i$. By Proposition 3.1, the stratifications $\mathcal{S}_k$ are
asymptotically holomorphic.
Applying Theorem 1.1 to these stratifications and starting from the
zero sections, we obtain asymptotically holomorphic sections of $E_k$ which 
are uniformly transverse to $\mathcal{S}_k$. The determinantal locus
$\Sigma_i(s_k)$ is precisely the set of points where the graph of $s_k$
intersects the stratum $S^i_k$. The result of uniqueness up to isotopy
is obtained by applying Theorem 3.2.
\end{proof}

However, our main application is that of maps to projective spaces.
Observe that, given a section $s=(s_1,\dots,s_{m+1})$ of a vector bundle of 
the form $\C^{m+1}\otimes L$, where $L$ is a line bundle over $X$, we can 
construct away from its zero set a projective map
$\mathbb{P}s=(s_1\!:\!\dots\!:\!s_{m+1}):X-s^{-1}(0)\to\CP^m$.

Recall that the space of jets of holomorphic maps from
$\C^n$ to $\C^m$ carries a natural partition into submanifolds,
the Boardman ``stratification'' \cite{arnold,boardman}. Restricting oneself 
to {\it generic} $r$-jets, the strata $\Sigma_I$, labelled by $r$-tuples 
$I=(i_1,\dots,i_r)$ 
with $i_1\ge\dots\ge i_r\ge 0$, are defined in the following way. 
Given a generic holomorphic map $f$, call
$\Sigma_i(f)$ the set of points where $\dim \mathrm{Ker}\, df=i$, and
denote by $\Sigma_i$ the set of holomorphic $1$-jets corresponding to such
points (i.e., $\Sigma_i$ is the set of $1$-jets $(\sigma_0,\sigma_1)$
such that $\dim\mathrm{Ker}\,\sigma_1=i$). The submanifolds $\Sigma_i$
determine a stratification of $\J^1_{n,m}$. For
a generic holomorphic map $f$ the critical loci $\Sigma_I(f)$ are 
smooth submanifolds defining a partition of $\C^n$.
Therefore, we can define 
inductively $\Sigma_{i_1,\dots,i_r}(f)$ as the set of points of
$\Sigma_{i_1,\dots,i_{r-1}}(f)$ where the kernel of the
restriction of $df$ to $T\Sigma_{i_1,\dots,i_{r-1}}(f)$ has dimension $i_r$
(in particular, $\Sigma_{i_1,\dots,i_{r-1},0}(f)$ is open in
$\Sigma_{i_1,\dots,i_{r-1}}(f)$ and corresponds to the set of points where
$f$ restricts to $\Sigma_{i_1,\dots,i_{r-1}}(f)$ as an immersion).

It is easy to check that the $r$-jet of $f$ at a given point of $\C^n$
completely determines in which $\Sigma_I(f)$ it lies~; therefore,
one can define $\Sigma_I\subset\J^r_{n,m}$ as the set of 
$r$-jets $j^r f(x)$ of generic holomorphic maps $f:\C^n\to\C^m$ at points 
$x\in\Sigma_I(f)$. In other terms, $\Sigma_I(f)=\{x\in\C^n,\ j^r f(x)\in
\Sigma_I\}$. It is a classical result \cite{boardman} that the $\Sigma_I$'s 
are  smooth submanifolds and define a partition of the space of generic 
holomorphic $r$-jets (an open subset in $\J^r_{n,m}$ whose complement has 
codimension $\ge n+1$), which can be extended into a partition of
$\J^r_{n,m}$ by smooth submanifolds.

The Boardman classes $\Sigma_I$ play a fundamental role in singularity
theory, and they completely determine the classification of 
singularities in certain dimensions. 
For low enough values of $r$, $m$ or $n$, the submanifolds $\Sigma_I$ define
a genuine stratification of the jet space $\J^r_{n,m}$. However, as observed 
by Boardman, things become more complicated as the dimension
increases, and the boundary of $\Sigma_I$ is in general not a union of entire 
strata~; in high dimensions Boardman classes do not even define a
quasi-stratification. 

Still, there exist well-known methods that allow
Boardman's partitions to be refined into finite Whitney stratifications of
$\J^r_{n,m}$. An example of such a construction can be found in the
work of Mather \cite{mather} (the constructed object
is tautologically a finite Whitney stratification, and one easy checks that
each Boardman class is a union of several of its strata).
\medskip

We now consider the case of maps to projective spaces defined by
asymptotically holomorphic sections of $E_k=\C^{m+1}\otimes L_k$ over $X$.
We want to construct a natural approximately holomorphic analogue of the 
Thom-Boardman stratifications, by defining certain submanifolds in $\J^r E_k$. 
In order to make things easier by avoiding a lengthy analysis of the
boundary structure at the points where the vanishing of the section prevents
the definition of a projective map, our aim will only be to construct
quasi-stratifications of $\J^r E_k$ rather than genuine stratifications.

We first define $Z=\{(\sigma_0,\dots,\sigma_r)\in\J^r E_k,\ \sigma_0=0\}$,
i.e.\ $Z$ is the set of $r$-jets of sections which vanish at the considered
point. As observed in \S 3.2, $\Theta_Z$ consists of all points of $Z$ such
that $\sigma_1$ is surjective. Next, observe that any point
$(\sigma_0,\dots,\sigma_r)\in\J^r E_k$ which does not belong to $Z$
determines the (symmetric) holomorphic $r$-jet $(\phi_0,\dots,\phi_r)$ of 
a map to $\CP^m$~: $\phi_0\in\CP^m$, $\phi_1\in T_x^*X^{1,0}\otimes
T_{\phi_0}\CP^m$, $\dots$, $\phi_r\in (T_x^*X^{1,0})^{\otimes
r}_\mathrm{sym}\otimes T_{\phi_0}\CP^m$ are defined in terms of
$\sigma_0,\dots,\sigma_r$ by expressions involving the projection map from
$\C^{m+1}-\{0\}$ to $\CP^m$ and its derivatives. In fact, one easily checks
that, if $(\sigma_0,\dots,\sigma_r)=j^r s$
is the symmetric holomorphic part of the $r$-jet of a section of $E_k$,
then $(\phi_0,\dots,\phi_r)=j^r f$ is the symmetric holomorphic part of the
$r$-jet of the corresponding projective map.
Using this notation, define 
$$\Sigma_i=\{(\sigma_0,\dots,\sigma_r)\in\J^r E_k,\ \sigma_0\neq 0,
\ \dim\mathrm{Ker}\,\phi_1=i\}.$$
For $\max(0,n-m)<i\le n$, one easily checks that $\Sigma_i$ is a smooth 
submanifold of $\J^rE_k$, and that $\partial\Sigma_i$ is the union of 
$\bigcup_{j>i} \Sigma_j$ 
and a subset of $Z-\Theta_Z$~: indeed, observe that if $n\ge m$, then 
for any $(\sigma_0,\dots,\sigma_r)\in\overline{\Sigma}_i\cap Z$ we have
$\dim\mathrm{Ker}\,\sigma_1\ge i-1>n-(m+1)$ and therefore $\sigma_1$
is not surjective, while in the case $n<m$ dimensional reasons prevent
$\sigma_1$ from being surjective.

Next, we assume that $r\ge 2$, and observe that $\Theta_{\Sigma_i}$ is
the set of points $(\sigma_0,\dots,\sigma_r)\in\Sigma_i$ such that 
$$\Xi_{i;(\sigma_0,\dots,\sigma_r)}=
\{u\in T_xX^{1,0},\ (\iota_u\sigma_1,\dots,\iota_u\sigma_r,0)\in
T_{(\sigma_0,\dots,\sigma_r)}\Sigma_i\}$$ has the expected codimension in
$T_xX^{1,0}$ (i.e., the same codimension as $\Sigma_i$ in $\J^r E_k$).
Indeed, by definition $(\sigma_0,\dots,\sigma_r)$ belongs to 
$\Theta_{\Sigma_i}$ if and only if the $(r+1)$-jet 
$(\sigma_0,\dots,\sigma_r,0)$, viewed as a $1$-jet in $\J^r E_k$, intersects
$\Sigma_i$ transversly (because the definition of 
$\Sigma_i$ involves only $\sigma_0$ and $\sigma_1$, the choice of a lift 
in $\J^{r+1}E_k$ does not matter, so we can choose the $(r+1)$-tensor 
component to be zero). By convention (see \S 3.2), this element of 
$\J^{r+1}E_k$ corresponds to the $1$-jet of a local section $\sigma$ of 
$\J^r E_k$ satisfying, at the given point $x\in X$, 
$\sigma(x)=(\sigma_0,\dots,\sigma_r)$ and
$\nabla\sigma(x)=(\sigma_1,\dots,\sigma_{r+1})$~: the covariant derivative
contains no antiholomorphic or antisymmetric terms. The graph of $\sigma$
intersects $\Sigma_i$ transversely if and only if $\{u\in TX,\
\nabla\sigma(x).u\in T_{\sigma(x)}\Sigma_i\}$ has the expected dimension,
hence the above criterion.

With this understood, we can define inductively, for $p+1\le r$,
$$\Sigma_{i_1,\dots,i_{p+1}}=\{\sigma\in \Theta_{\Sigma_{i_1,\dots,i_p}},
\ \dim(\mathrm{Ker}\,\phi_1\cap \Xi_{(i_1,\dots,i_p);\sigma})=i_{p+1}\},$$
where $\Xi_{I;\sigma}=\{u\in T_xX^{1,0},
\ (\iota_u\sigma_1,\dots,\iota_u\sigma_r,0)\in
T_{(\sigma_0,\dots,\sigma_r)}\Sigma_I\}$
as above, and $\Theta_{\Sigma_I}$ again consists of all points 
$\sigma\in\Sigma_I$ such that $\Xi_{I;\sigma}$ has the same codimension
in $T_xX^{1,0}$ as $\Sigma_I$ in $\J^r E_k$.

For $i_1\ge\dots\ge i_{p+1}\ge 1$,
$\Sigma_{i_1,\dots,i_{p+1}}$ is a smooth submanifold in $\J^r
E_k$, and its closure inside $\Sigma_{i_1,\dots,i_p}$ is obtained
by adding $\bigcup_{j>i_{p+1}}\Sigma_{i_1,\dots,i_p,j}$ and a subset of
$\Sigma_{i_1,\dots,i_p}-\Theta_{\Sigma_{i_1,\dots,i_p}}$. However, it is
quite difficult to fully understand the boundary structure of 
$\Sigma_{i_1,\dots,i_{p+1}}$~; the situation is exactly
the same as in standard Boardman theory for holomorphic jets, except that,
besides pieces of $\Sigma_{j_1,\dots,j_q}$ where $q\le p+1$
and $(j_1,\dots,j_q)\ge (i_1,\dots,i_q)$ for the lexicographic order, 
the boundary of $\Sigma_{i_1,\dots,i_{p+1}}$ also contains a subset of
$Z-\Theta_Z$. 

In low dimensions and/or for low values of $r$, it can be checked that the 
submanifolds $Z,\Sigma_i,\Sigma_{i_1,i_2},\dots,\Sigma_{i_1,\dots,i_r}$ 
determine a finite Whitney quasi-stratification of $\J^r E_k$~; for example
when $r=1$ this is an immediate consequence of the above discussion.

However, in larger dimensions it is necessary to refine Boardman's
construction as in the holomorphic case. The important observation is that,
when $\J^r E_k$ is trivialized by choosing local asymptotically holomorphic 
coordinates and sections, the partition of $\J^r E_k-Z$ described above
corresponds exactly to the partition of the space of $r$-jets of maps to 
$\CP^m$ given by Boardman classes. Therefore, we can circumvent the problem
by refining Boardman's partition of $\J^r_{n,m}$ into a genuine 
stratification as explained above, lifting it by the projectivization map to
a stratification of the space of non-vanishing jets in $\J^r_{n,m+1}$, and
finally pull it back to obtain a stratification of $\J^r E_k-Z$. As in the
holomorphic case, the $\Sigma_I$ classes are realized as unions of strata~;
therefore, transversality to this stratification implies transversality to
the $\Sigma_I$'s. Moreover, all strata (except for the open one which we 
discard anyway) are contained in the closure of $\Sigma_1$, so that
the boundary structures near $Z$ are entirely contained in $Z-\Theta_Z$~;
therefore adding $Z$ to this stratification yields a
quasi-stratification of $\J^r E_k$.

\begin{defi}
Given asymptotically very ample line bundles $L_k$ over 
the manifold
$(X^{2n},J)$, and 
setting $E_k=\C^{m+1}\otimes L_k$, the {\em Boardman stratification} of 
$\J^r E_k$ is the quasi-stratification given by the submanifold $Z$ and
by a refined Thom-Boardman stratification of $\J^r E_k-Z$.
\end{defi}

\begin{coro}
Let $(X,\omega)$ be a compact symplectic manifold endowed with an $\omega$-tame
almost-complex structure $J$, let $L_k$ be an asymptotically very ample 
sequence of line bundles over $(X,J)$, and let $E_k=\C^{m+1}\otimes L_k$.
Then, for all large enough values of $k$ there exist asymptotically 
holomorphic sections $s_k$ of $E_k$ such that the $r$-jets $j^r s_k$ are
uniformly transverse to the Boardman stratifications of $\J^r E_k$.

In particular, the zero sets $Z_k=s_k^{-1}(0)$ are smooth symplectic
codimension $2m$ submanifolds in $X$, and the holomorphic $r$-jets of
the projective maps $f_k=\mathbb{P}s_k:X-Z_k\to\CP^m$ behave at every point
in a manner similar to those of generic holomorphic maps from a complex
$n$-fold to $\CP^m$. Moreover, the singular loci $\Sigma_I(f_k)=
\{x\in X-Z_k,\ j^r f_k(x)\in \Sigma_I\}$ are smooth symplectic submanifolds
of the expected codimension and define a partition of $X-Z_k$.
Finally, the sections $s_k$ and the maps $f_k$
are, for large $k$, canonical up to isotopy, independently of the chosen
almost-complex structure on~$X$.
\end{coro}%

\begin{proof} By construction the Boardman stratifications of $\J^r E_k$
satisfy the assumptions of Proposition 3.1, as in every fiber of $\J^r E_k$
they can be identified with the same holomorphic quasi-stratification of 
$\J^r_{n,m+1}$.
As a consequence, they are asymptotically holomorphic, and the existence 
of asymptotically holomorphic sections of $E_k$ with the desired 
transversality properties is an immediate consequence of Theorem 1.1. 
The properties of $Z_k$ follow immediately from the uniform transversality 
to the stratum $Z$ of vanishing sections, while the properties of $f_k$ 
are direct consequences of the uniform transversality to the Boardman
strata (recall that each $\Sigma_I$ is smooth and is a union of strata). 
Finally, the uniqueness result is obtained by applying Theorem 3.2.
\end{proof}

Corollary 5.3 is, in a certain sense, a fundamental result of asymptotically
holomorphic singularity theory. Still, it falls short of the natural goal
that one may have in mind at this point, namely the construction of
approximately holomorphic projective maps which are near every point of $X$
topologically conjugate in approximately holomorphic coordinates to generic
holomorphic maps between complex manifolds. 

Indeed, in order to achieve such a result, one needs to obtain some control on 
the antiholomorphic part of the jet of $f_k$ at the points of the singular
loci $\Sigma_I(f_k)$~: roughly speaking, $\dbar f_k$ must be much smaller
than $\partial f_k$ in every direction and at every point, and when
$\partial f_k$ is singular this is no longer an immediate consequence of
asymptotic holomorphicity and transversality. Note however that the behavior 
of $f_k$ near the set of base points $Z_k$ is always the expected one.

In many cases, it is possible to perturb slightly the sections $s_k$ (by
less than a fixed multiple of $c_k^{-1/2}$, which affects neither
holomorphicity nor transversality properties) along the singular loci
in order to obtain the proper topological picture for $f_k$.

The easiest case is $m\ge 2n$, where it is enough to consider $1$-jets,
and all the strata turn out to be of codimension greater than $n$~; the 
uniform transversality of $s_k$ to the Boardman stratification then implies 
that the maps $f_k$ are approximately holomorphic immersions. 
Moreover, when $m\ge 2n+1$ an arbitrarily small perturbation is 
enough to get rid of multiple points, thus giving approximately holomorphic
embeddings into projective spaces, a result already obtained by Mu\~noz, 
Presas and Sols \cite{mps}.

Next, we can consider the case $m=1$, where $1$-jets are again sufficient,
and the only interesting Boardman stratum is $\Sigma_n$, of complex
codimension $n$, corresponding to critical points of $\CP^1$-valued maps.
The sections of $\C^2\otimes L_k$ given by Corollary 5.3 vanish along 
smooth codimension 4 base loci~; moreover,
the differential $\partial f_k$ of the $\CP^1$-valued map $f_k$ only
vanishes at isolated points, and does so in a non-degenerate way. These
transversality properties are precisely those imposed by Donaldson in his
construction of symplectic Lefschetz pencils \cite{D2}; the only missing
ingredient is an extra perturbation near the zeroes of $\partial f_k$ in 
order to get rid of the antiholomorphic terms and therefore ensure that
they are genuine non-degenerate critical points, thus making $f_k$ a complex
Morse function.

The last case we will consider is when $m=2$. In this case, we need to
consider $2$-jets, and the relevant Boardman strata are $\Sigma_{n-1}$,
of complex codimension $n-1$, and $\Sigma_{n-1,1}$, of complex codimension
$n$ (the other strata have codimension greater than $n$). 
The sections of $\C^3\otimes L_k$ constructed by Corollary 5.3 vanish along
smooth codimension 6 base loci. The $\CP^2$-valued maps $f_k$ are
submersions outside of the smooth symplectic curves
$R_k=\Sigma_{n-1}(f_k)$, and
the restriction of $f_k$ to $R_k$ is an immersion except at the points of
$C_k=\Sigma_{n-1,1}(f_k)$. After a suitable perturbation in 
order to ensure the vanishing of some antiholomorphic derivatives of $f_k$ 
along $R_k$, one obtains a situation similar to that described in previous
papers \cite{A2,A3}: at every point of $R_k-C_k$, a local model
for $f_k$ in approximately holomorphic coordinates is
$(z_1,\dots,z_n)\mapsto (z_1^2+\dots+z_{n-1}^2,z_n)$, while at the points
of $C_k$ the local model becomes $(z_1,\dots,z_n)\mapsto
(z_1^3+z_1z_n+z_2^2+\dots+z_{n-1}^2,z_n)$ and the symplectic curve
$f_k(R_k)\subset\CP^2$ presents an isolated cusp singularity.

In the general case, the most promising strategy to achieve topological 
conjugacy to generic holomorphic local models is to perturb the
sections $s_k$ in order to make sure that, along each stratum $\Sigma_I(f_k)$,
the germ of $f_k$ is holomorphic along the normal directions to 
$\Sigma_I(f_k)$. Such perturbations should be relatively easy to construct by 
methods similar to those in the above-mentioned papers\cite{A2,A3}, provided
that one starts from the strata of lowest dimension. This approach will
be developped in a forthcoming paper.

Finally, let us formulate some natural extensions of Corollary 5.3 to more
general situations. First, we mention the case when the asymptotically
very ample line bundles $L_k$ are replaced by vector bundles of rank 
$\nu\ge 2$. In that case, and provided that $m\ge\nu$, the projective maps 
defined by sections $s_k$ of $\C^{m+1}\otimes L_k$ are replaced by maps 
$\mathrm{Gr}(s_k)$ taking values in the Grassmannian $\mathrm{Gr}(\nu,m+1)$ 
of $\nu$-planes in $\C^{m+1}$, defined at every point of $X$ where the $m+1$ 
chosen sections generate the whole fiber of $L_k$. More precisely, at every
such point there exist $m+1-\nu$ independent linear relations between the 
$m+1$ components $s_k^1,\dots,s_k^{m+1}$, and these $m+1-\nu$ linear
equations in $m+1$ variables determine a $\nu$-dimensional complex subspace
$\mathrm{Gr}(s_k)$ in $\C^{m+1}$. By adapting Corollary 5.3 to this
situation, it is for example possible to recover the Grassmannian embedding
result of Mu\~noz, Presas and Sols \cite{mps}.

Another direction in which Corollary 5.3 can be improved is by adding extra
transversality requirements to the projective maps $f_k$. For example,
given a stratified holomorphic submanifold $\mathcal{D}=(D_a)_{a\in A}$ in 
$\CP^m$, we can require the transversality of the map $f_k$ to
$\mathcal{D}$. Indeed, $\mathcal{D}$ induces a stratification
$\tilde{\mathcal{D}}_k$ of $\J^r E_k$, in which each stratum consists of the 
jets $(\sigma_0,\dots,\sigma_r)$ such that $\mathbb{P}\sigma_0$ belongs
to a certain stratum $D_a$ of $\mathcal{D}$ (in fact, this 
stratification only involves the $0$-jet part). Starting from sections
$s_k$ of $E_k$ given by Corollary 5.3, we can apply Theorem 1.1 to the
stratifications $\tilde{\mathcal{D}}_k$ (which one easily shows to 
to be asymptotically holomorphic by Proposition 3.1)~; this yields
asymptotically holomorphic sections $\tilde{s}_k$ which are uniformly 
transverse to $\tilde{\mathcal{D}}_k$ but differ from $s_k$ by an amount small
enough to ensure that the transversality of the jets to the Boardman
stratification is preserved. In this way, one obtains projective
maps which have the same properties as in Corollary 5.3 and additionally
are uniformly transverse to the stratified submanifold $\mathcal{D}$.
This extends a result of Mu\~noz, Presas and Sols \cite{mps} where
asymptotically holomorphic embeddings are made transverse to a given
submanifold of $\CP^m$.

Another class of stratifications of $\J^rE_k$ that we can consider are 
those obtained from lower-dimensional Boardman stratifications by linear
projections. Namely, fix $q<m$, and let $\pi:\C^{m+1}\to\C^{q+1}$ be a
linear projection~; $\pi$ induces maps $\tilde\pi:\J^r(\C^{m+1}\otimes L_k)
\to\J^r(\C^{q+1}\otimes L_k)$, and the inverse images by $\tilde\pi$ of the
Boardman stratifications of $\J^r(\C^{q+1}\otimes L_k)$ are asymptotically
holomorphic quasi-stratifications of $\J^r(\C^{m+1}\otimes L_k)$. The
transversality of $j^r s_k$ to these quasi-stratifications is equivalent to 
that of $j^r(\pi(s_k))$ to the Boardman stratifications~; 
denoting by $\bar\pi$ the map from $\CP^m$ to $\CP^q$ induced by $\pi$, this
is also equivalent to the genericity of the holomorphic jets of the
projective maps $\bar\pi\circ f_k$. Therefore, by applying Theorem 1.1 as in 
the previous example, we can obtain projective maps $f_k$ with the same 
genericity properties as in Corollary 5.3 and such that the
maps $\bar\pi\circ f_k$ also enjoy similar properties. Even
better, by iteratedly applying Theorem 1.1 we can obtain the same property
for any given finite family of linear projections. For example, when $m=2$ and
considering projections of $\C^3$ to $\C^2$ along coordinate axes, one 
obtains exactly the transversality properties which are needed in order to 
extend Moishezon-Teicher braid group techniques to the study of symplectic 
manifolds \cite{AK,A3}.

To conclude, let us mention a different class of potential applications
of Theorem 1.1, following the ideas of Donaldson and Smith. 
As shown by Donaldson \cite{D2}, any compact
symplectic 4-manifold carries structures of symplectic Lefschetz pencils
obtained from pairs of sections of asymptotically very ample line bundles
$L_k$~; after blowing up the base points, we obtain Lefschetz fibrations 
over $\CP^1$, which may also be thought of as maps from $\CP^1$ to the 
moduli space $\bar{M}_g$ of stable curves of a certain genus $g$. These 
maps become asymptotically holomorphic as one considers pencils given by 
sections of $L_k$ for $k\to+\infty$. In a largely unexplored class of
constructions, one considers certain vector bundles over $\CP^1$ naturally
arising from the Lefschetz fibrations~: for example, spaces of holomorphic
sections of certain bundles over each fiber,
or pull-backs by the maps from $\CP^1$ to $\bar{M}_g$ of vector
bundles over $\bar{M}_g$. It often turns out that these bundles over $\CP^1$
either are naturally asymptotically very ample or become so after tensor
product by the line bundles $O(k)$. Theorem 1.1 can then used in order to
obtain sections with suitable genericity properties, which in turn give
rise to interesting geometric or topological structures. In some cases
the objects naturally arising are sheaves rather than bundles, but the
same type of argument should remain valid. It is to be
expected that some interesting results about symplectic 4-manifolds and
Lefschetz pencils can be obtained in this way, as similar considerations
(but at a much more sophisticated level) have for example led to Donaldson
and Smith's proof of the existence of a pseudo-holomorphic curve realizing
the canonical class via Lefschetz fibrations \cite{ds}.

\end{document}